\documentclass[12pt]{amsart}
\usepackage{amsfonts, amssymb, latexsym, epsfig}

\usepackage{amsthm}
\usepackage{amsmath}
\usepackage{setspace}
\usepackage{graphicx}
\usepackage[usenames, dvipsnames]{xcolor}
\usepackage{tikz}
\usetikzlibrary{calc, shapes, backgrounds,arrows,positioning}
\tikzset{>=stealth',
  head/.style = {fill = white, text=black},
  plaque/.style = {draw, rectangle, minimum size = 10mm}, 
  pil/.style={->,thick},
  junct/.style = {draw,circle,inner sep=0.5pt,outer sep=0pt, fill=black}
  }
 \definecolor{ttffcc}{rgb}{0.2,1.,0.8}
\definecolor{aqaqaq}{rgb}{0.6274509803921569,0.6274509803921569,0.6274509803921569}
\definecolor{uuuuuu}{rgb}{0.26666666666666666,0.26666666666666666,0.26666666666666666}
\definecolor{ffffff}{rgb}{1.,1.,1.}
\usepackage{enumerate}
\usepackage{ytableau}
\usepackage{hyperref}

\newtheorem{theorem}{Theorem}[section]
\newtheorem{lemma}[theorem]{Lemma}
\newtheorem{conjecture}[theorem]{Conjecture}
\newtheorem{corollary}[theorem]{Corollary}
\newtheorem{proposition}[theorem]{Proposition}

\theoremstyle{definition}

\newtheorem{remark}[theorem]{Remark}
\newtheorem{example}[theorem]{Example}

\theoremstyle{remark}

\makeatletter
\@addtoreset{case}{section}
\makeatother

\numberwithin{subcase}{case}

\numberwithin{subsubcase}{subcase}

\newcommand{\inc}[1]{\mathtt{Inc}(#1)}
\newcommand{\incr}[2]{\mathtt{IncRect}_#2(#1)}

\newcommand{\blank}{\phantom{2}}

\newcommand{\genomictableau}[2]{\mathtt{Gen}_{#2}(#1)}
\newcommand{\Pgenomictableau}[2]{\mathtt{PGen}_{#2}(#1)}
\newcommand{\Qgenomictableau}[2]{\mathtt{QGen}_{#2}(#1)}
\newcommand{\ballot}[2]{\mathtt{Ballot}_{#2}(#1)}
\newcommand{\Pballot}[2]{\mathtt{PBallot}_{#2}(#1)}
\newcommand{\Qballot}[2]{\mathtt{QBallot}_{#2}(#1)}
\newcommand{\pieri}[2]{\mathtt{PF}_#2(#1)}
\newcommand{\jdt}[2]{\mathtt{jdt}_#2(#1)}

\newcommand{\im}{\mathrm{im \;}}
\newcommand*{\Scale}[2][4]{\scalebox{#1}{$#2$}}%
\newcommand*\circled[1]{\tikz[baseline=(char.base)]{
  \node[shape=circle,draw,inner sep=1pt] (char) {$#1$};}}

\newcommand{\grey}{\textcolor{gray}}

\newcommand{\gap}{\hspace{1in} \\ \vspace{-.2in}}

\newcommand{\Z}{\mathbb{Z}}

\newcommand{\GG}{\mathcal{G}}

\newcommand{\excise}[1]{}

\begin{document}
\pagestyle{plain}

\title{Genomic tableaux}
\author{Oliver Pechenik}
\author{Alexander Yong}
\address{Department of Mathematics \\ University of Illinois at Urbana--Champaign \\ Urbana, IL 61801 \\ USA}
\email{pecheni2@illinois.edu, ayong@uiuc.edu}
\date{March 28, 2016}
\vspace{-.1in}
\begin{abstract}
We explain how \emph{genomic tableaux} [Pechenik-Yong '15] are a semistandard complement to \emph{increasing tableaux} [Thomas-Yong '09]. 
From this perspective, one inherits genomic versions of \emph{jeu de taquin}, Knuth equivalence,  infusion and 
Bender-Knuth involutions, as well as Schur functions from (shifted) semistandard Young tableaux theory. 
These are applied to obtain
new Littlewood-Richardson rules for  $K$-theory Schubert calculus of Grassmannians (after [Buch '02]) and 
maximal orthogonal Grassmannians (after [Clifford-Thomas-Yong '14], [Buch-Ravikumar '12]). For the unsolved case of 
Lagrangian Grassmannians, sharp upper and lower bounds using genomic tableaux are conjectured. 
\end{abstract}
\maketitle 

\vspace{-.2in}
\tableofcontents
\vspace{-.4in}

\ytableausetup{boxsize=1.1em}

\section{Introduction}

\subsection{History and overview}
Let ${\sf Sym}$ be the ring of symmetric functions. Textbook theory 
of ${\sf Sym}$ concerns the basis of Schur functions and the package of Young tableau
algorithms for which the Littlewood-Richardson rule
is a centerpiece. Interpreting these
polynomials in the Schubert calculus of Grassmannians one is led, via a $K$-theoretic
generalization, to \emph{symmetric Grothendieck 
functions} $\{G_{\lambda}\}$, a deformation of the Schur basis $\{s_{\lambda}\}$. This line of inquiry started with 
\cite{anneaux}. The first combinatorial rule for $G_{\lambda}$
was given by \cite{Fomin.Kirillov} whereas the first tableau formula was
found by \cite{Buch:KLR}.

There is interest in finding $K$-analogues of elements of the classical Young tableau theory; see, e.g.,
\cite{Lenart, Buch:KLR, BKSTY, Thomas.Yong:V, Buch.Samuel, MinnREU, PaPyII, MinnREU2, Morse}. Although the Grothendieck functions were originally studied for 
geometric reasons, the  
combinatorics has been part of a broader conversation in algebraic and enumerative combinatorics, e.g., \emph{Hopf algebras} \cite{LaPy07, PaPy, Patrias},
\emph{cyclic sieving} \cite{Pe:cyclic, Rhoades, Stokke},
\emph{Demazure characters} \cite{Monical},
\emph{homomesy} \cite{BPS},
longest increasing subsequences of random words \cite{TY:LIS},
poset edge densities \cite{Reiner.Tenner.Yong}, and 
\emph{plane partitions} \cite{DPS, HPPW}.

In \cite{Thomas.Yong:V}, a \emph{jeu de taquin} theory for
\emph{increasing tableaux} was introduced. These tableaux are fillings of 
Young diagrams $\nu / \lambda$ with $[\ell]:=1,2,\ldots,\ell$ where $\ell\leq |\nu / \lambda|$ and the entries increase in rows and columns (labels may be repeated). If $\ell=|\nu / \lambda|$, these are standard Young tableaux and 
increasing tableau results closely parallel those for
standard Young tableaux. An outcome was a new Littlewood-Richardson rule for $\{G_{\lambda}\}$
(after \cite{Buch:KLR}) and its 
\emph{minuscule} extension (see \cite{Thomas.Yong:adv,
Buch.Ravikumar, Clifford.Thomas.Yong, Buch.Samuel}).

In \cite{Thomas.Yong:H_T}, a \emph{jeu de taquin}-based  
Littlewood-Richardson rule for torus-equivariant $K$-theory of Grassmannians was conjectured. 
In \cite{PY:full}, we proved this conjecture by defining \emph{genomic tableaux} as 
a semistandard analogue of increasing tableaux.

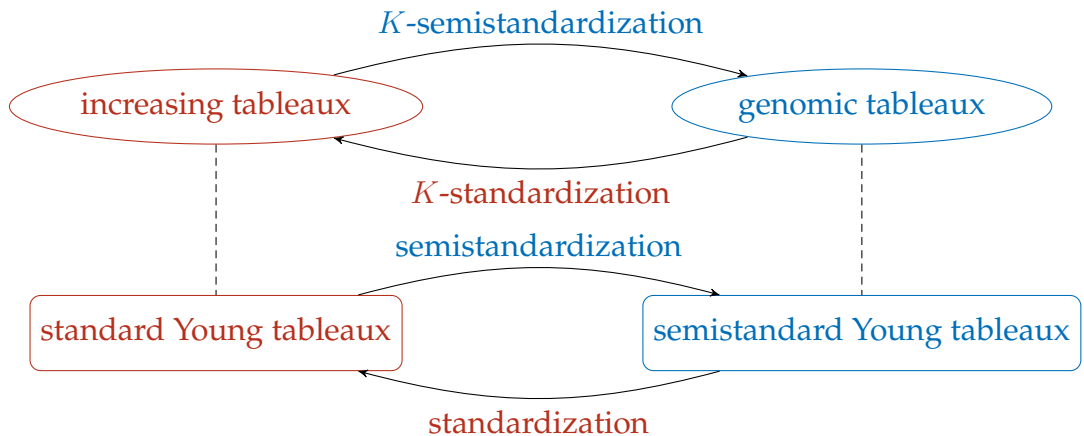
\vspace{-.2in}
\[\begin{tikzpicture}[node distance=1cm, auto]
 \node[plaque, BrickRed, ellipse] (A) {increasing tableaux};
 \node[right=3.3cm of A, plaque, NavyBlue, ellipse] (B) {genomic tableaux};
 \node[below=2.0cm of A, plaque, rounded corners, BrickRed] (C) {standard Young tableaux};
 \node[below=2.0cm of B, plaque, rounded corners, NavyBlue] (D) {semistandard Young tableaux};
  \begin{scope}[nodes = {draw = none}]
    \path (A) edge[bend left=15, ->] node[midway, above, NavyBlue]  {$K$-semistandardization}  (B)
      (A) edge[bend right=15, <-] node[midway,below, BrickRed] {$K$-standardization} (B)
      (A) edge [densely dashed] (C)
      (C) edge[bend left=15, ->] node[midway, above, NavyBlue]  {semistandardization} (D)
      (C) edge[bend right=15, <-] node[midway,below, BrickRed] {standardization} (D)
      (B) edge [densely dashed] (D)
      ;
  \end{scope}
\end{tikzpicture}
\]

Our goal is a theory of genomic tableaux parallel to that of 
\cite{Thomas.Yong:V} for increasing tableaux. The Schubert calculus
application in \cite{PY:full, PY:puzzles}
used \emph{edge-labeled} genomic tableaux. However, in anticipation of other applications, we give a logically independent development of 
genomic tableau combinatorics in the basic (i.e., non-edge labeled) case
and in the shifted setting. 
The first applications are to give
new Littlewood-Richardson-type rules for (ordinary) $K$-theory of Grassmannians and maximal orthogonal Grassmannians. Furthermore, modifications of these rules give conjectural upper and lower bounds for the $K$-theory 
structure constants of Lagrangian Grassmannians. 

\subsection{Genomic tableau results}
\label{sec:preview}
Let $S$ be a semistandard Young tableau of a shape $\nu/\lambda$.  Place a total order on those boxes with entry $i$ using left to right order.
A {\bf gene} $\GG$ (of family $i$) is a collection of consecutive boxes in this order, where no two lie in the same row; we write ${\tt family}(\GG)=i$.
A {\bf genomic tableau} $T$ is a semistandard tableau 
together with a partition of its boxes into genes. We
indicate the partition by color-coding the boxes.
The {\bf content} of $T$ is the number of genes of each family. 
Note, a semistandard
tableau $T$ is a genomic tableau where each gene is a single box.
Moreover, the content of $T$ agrees with the usual notion for semistandard tableaux.

\begin{example}\label{ex:typeA_genomic tableau}
$T =
\begin{ytableau}
 *(lightgray)\blank & *(lightgray)\blank &*(red) 1 &*(green) 2\\
*(SkyBlue) 1 & *(red) 1 &*(green) 2 \\
*(green) 2
\end{ytableau}$ \ has content $(2, 1)$ since there are two genes of family $1$ and one of family $2$.\qed
\end{example}

A {\bf genotype} $G$ of a genomic tableau $T$ is a choice of a single box from each gene.\footnote{The genomic analogy is that boxes of a gene are \emph{alleles} and the other genes of the same family are \emph{paralogs}.} We depict $G$ by erasing the entries in all unchosen boxes of $T$.

\begin{example}\label{ex:typeA_genomic tableauII}
Continuing Example~\ref{ex:typeA_genomic tableau},
\[\begin{ytableau}
 *(lightgray)\blank & *(lightgray)\blank & \blank & \blank \\
 1 & 1 & \blank \\
 2
\end{ytableau} 
\hspace{.3cm}
\begin{ytableau}
 *(lightgray)\blank & *(lightgray)\blank & 1 & \blank \\
 1 & \blank & \blank \\
 2
\end{ytableau}
\hspace{.3cm}
\begin{ytableau}
 *(lightgray)\blank & *(lightgray)\blank & \blank & \blank\\
 1 & 1 & 2 \\
\blank
\end{ytableau}
\hspace{.3cm}
\begin{ytableau}
 *(lightgray)\blank & *(lightgray)\blank & 1 & \blank\\
 1 & \blank & 2 \\
\blank
\end{ytableau}
\hspace{.3cm}
\begin{ytableau}
 *(lightgray)\blank & *(lightgray)\blank & \blank & 2\\
 1 & 1 & \blank \\
\blank
\end{ytableau}
\hspace{.3cm}
\begin{ytableau}
 *(lightgray)\blank & *(lightgray)\blank & 1 & 2\\
 1 & \blank & \blank \\
\blank
\end{ytableau}\]
are the six genotypes of $T$.\qed
\end{example}

Suppose $U$ is any filling of a subset of boxes of a 
shape. The {\bf sequence}
${\tt seq}(U)$ of $U$ is the reading word obtained by 
reading its entries along rows from right to left and from top to bottom (ignoring empty boxes). Now, ${\tt seq}(U)$ is a {\bf ballot sequence}
if the number of $i$'s that appear is always weakly greater than the number of $(i+1)$'s that appear, at any point in the sequence. A genomic tableau $T$ is {\bf ballot} if ${\tt seq}(G)$ is a ballot sequence for every genotype $G$ of $T$. Notice if each gene of $T$ is a single box, there is a unique genotype (namely, the underlying semistandard
tableau of $T$) and the concept of a ballot tableau coincides with the same notion for semistandard tableaux.

\begin{example}
The genotypes of
Example~\ref{ex:typeA_genomic tableauII}
respectively have sequences: $112$, $112$, $211$, $121$, $211$, and $211$. Since $211$ is not a ballot sequence, $T$ is not ballot. \qed
\end{example}

Our results are:
\begin{enumerate} 
\item A $K$-analogue of the \emph{(semi)standardization maps} between standard and semistandard tableaux. This 
relates genomic tableaux to increasing tableaux. 
\item Using (1), we acquire genomic analogues of \emph{Knuth equivalence}, \emph{jeu de taquin}, \emph{infusion} and 
\emph{Bender-Knuth involutions}. 
\item Using (2), we describe a new
basis $\{U_{\lambda}\}$ of ${\sf Sym}$ where each $U_{\lambda}$ is a generating series over genomic
tableaux of shape $\lambda$. This is a deformation of the Schur basis.
\item We give shifted analogues of (1)--(3).
\end{enumerate}

\subsection{Genomic rules in Schubert calculus} 
\label{sec:result_A}
Let 
\[X={\rm Gr}_k({\mathbb C}^n)\] 
denote the {\bf Grassmannian} of $k$-dimensional linear subspaces of ${\mathbb C}^n$. The
general linear group ${\sf GL}_n$, and its Borel subgroup ${\sf B}_{-}$ of lower triangular invertible matrices, act on $X$ by change of basis. 
This action decomposes $X$ into ${\sf B}_{-}$-orbits
\[X_{\lambda}^{\circ}\cong {\mathbb C}^{k(n-k)-|\lambda|}\] (the {\bf Schubert cells});
here $\lambda$ is a Young diagram contained in the rectangle $k\times (n-k)$. 
The 
Zariski closure of $X_{\lambda}^\circ$ is the 
{\bf Schubert variety} 
\[X_{\lambda}=\bigsqcup_{\mu\supseteq \lambda} X_{\mu}^{\circ}.\]

Textbook discussion of Schubert calculus revolves around classes of $X_{\lambda}$ in the cohomology ring $H^{\star}(X,{\mathbb Z})$; see, e.g., \cite{Fulton:YT}. These classes form a ${\mathbb Z}$-linear basis of $H^{\star}(X,{\mathbb Z})$.
Their structure constants 
\[[X_{\lambda}]\cup[X_{\mu}]=\sum_{\nu}C_{\lambda,\mu}^{\nu}[X_{\nu}]\]
with respect to the cup product are given by the classical \emph{Littlewood-Richardson rule} that governs
the multiplication of Schur functions. Geometrically,
\[C_{\lambda,\mu}^{\nu}=\#\text{points in $g_1\cdot X_{\lambda}\cap g_2 \cdot X_{\mu}\cap g_3 \cdot X_{\nu^{\vee}}$}\] 
for generic $g_1,g_2,g_3\in {\sf GL}_n$
when this number is finite; it is zero otherwise. Here,
\[\nu^{\vee}=(n-k-\nu_k,n-k-\nu_{k-1},\ldots,n-k-\nu_1)\] 
is the rotation of the complement of $\nu$ in $k\times (n-k)$. 

There has been significant attention on 
$K$-theoretic Schubert calculus, which provides a richer setting for study; see, e.g., \cite{Brion, Buch:combinatorialKtheory, Vakil, Knutson:14} and the references therein. 
Recall that the {\bf Grothendieck ring} $K^{0}(X)$ is the free abelian group
generated by isomorphism classes $[V]$ of algebraic vector bundles over $X$ under the relation 
\[[V]=[U]+[W]\] 
whenever there is a short exact sequence 
\[0\to U\to V\to W\to 0.\] 
The product structure on $K^{0}(X)$ is given by the tensor product of vector bundles, i.e.,
\[[U]\cdot [V]=[U\otimes V].\] 
Since $X$ is a smooth projective variety, the structure sheaf ${\mathcal O}_{X_{\lambda}}$ has a resolution 
\[0\to V_N\to V_{N-1}\to \cdots \to V_1\to V_0\to {\mathcal O}_{X_{\lambda}}\to 0\]
by locally free sheaves. Therefore it makes sense to define the class $[{\mathcal O}_{X_{\lambda}}]$ by 
\[[{\mathcal O}_{X_{\lambda}}]:=\sum_{j=0}^N(-1)^j[V_j]\in K^{0}(X).\] 

Now, $\{[{\mathcal O}_{X_\lambda}]\}$ forms a ${\mathbb Z}$-linear basis of $K^0(X)$. Thus, define structure constants by
\begin{equation}
\label{eqn:structure}
[{\mathcal O}_{X_{\lambda}}]\cdot [{\mathcal O}_{X_{\mu}}]=\sum_{\nu} a_{\lambda,\mu}^{\nu}[{\mathcal O}_{X_{\nu}}].
\end{equation}
A.~Buch \cite{Buch:KLR} gave a combinatorial rule for
$a_{\lambda,\mu}^{\nu}$, thereby establishing
\[(-1)^{|\nu|-|\lambda|-|\mu|} a_{\lambda,\mu}^{\nu}\geq 0.\] 
A number of other rules have been
discovered since, see, e.g., \cite{Vakil, BKSTY, Thomas.Yong:V} and the references therein, as well as the references above. 

\begin{theorem}[Genomic Littlewood-Richardson rule]\label{thm:A_lr_rule}
$a_{\lambda, \mu}^\nu=(-1)^{|\nu| - |\lambda| - |\mu|}$ times the number of ballot genomic tableaux of shape $\nu / \lambda$ and content $\mu$.
\end{theorem}

Actually, in the case $|\nu| = |\lambda| + |\mu|$, $a_{\lambda,\mu}^{\nu}=C_{\lambda,\mu}^{\nu}$ and Theorem~\ref{thm:A_lr_rule}
recovers the original rule of D.E.~Littlewood-A.R.~Richardson for multiplication of Schur functions \cite{Littlewood.Richardson}.

\begin{example}
The tableau \
$\ytableaushort{ {*(lightgray)\blank} {*(lightgray)\blank} {*(red) 1}, {*(lightgray)\blank} {*(red) 1}, {*(green) 2}}$ \
is the unique witness of $a_{(2,1),(1,1)}^{(3,2,1)}=-1$. 
\qed
\end{example}

Using the tableau results outlined in Section~\ref{sec:preview}, the
proof of Theorem~\ref{thm:A_lr_rule} is derived from the corresponding Littlewood-Richardson
rule in \cite{Thomas.Yong:V}. 

Similarly, using our theorems on shifted genomic tableaux, we obtain a combinatorial rule for
the $K$-theoretic structure constants for maximal orthogonal Grassmannians from the corresponding 
theorem in \cite{Clifford.Thomas.Yong} (originally conjectured in \cite{Thomas.Yong:V}).  

In cohomology, there is a simple relation between the structure constants for maximal orthogonal and Lagrangian Grassmannians.
Hence, given the aforementioned results, it may come as a bit of surprise that there is no known combinatorial rule 
for the $K$-theory structure constants of Lagrangian Grassmannians. 
We use genomic tableaux to contribute related conjectural upper and lower bounds for these numbers.
We further 
conjecture an inequality between the structure constants of Lagrangian Grassmannians and the maximal orthogonal
Grassmannians. 

\section{K-(semi)standardization maps}
\label{sec:K-maps}
Let 
\[\genomictableau{\nu / \lambda}{\mu}=\{
\text{genomic tableaux of shape $\nu / \lambda$ with content $\mu = (\mu_1, \mu_2 \dots, \mu_{\ell(\mu)})$}\},\] 
and
\[\inc{\nu / \lambda}=\{\text{increasing tableaux of shape $\nu / \lambda$}\}.\]

 Define an order on the genes of $T \in \genomictableau{\nu / \lambda}{\mu}$ by 
$\GG_1 < \GG_2$ if ${\tt family}(\GG_1) < {\tt family}(\GG_2)$ or if ${\tt family}(\GG_1) = {\tt family}(\GG_2)$ with all boxes of $\GG_1$ west of all boxes of $\GG_2$. 

\begin{lemma}\label{lem:gene_total_order}
The order $<$ on genes of $T$ is a total order.
\end{lemma}
\begin{proof}
When showing two genes $\GG_1$ and $\GG_2$ are comparable in the order $<$, the only concern is if ${\tt family}(\GG_1) = {\tt family}(\GG_2) = k$. By definition, a gene of family $k$ consists of boxes of entry $k$ that are consecutive in the left to right order on such boxes. Hence either all boxes of $\GG_1$ are west of the boxes of $\GG_2$ or \emph{vice versa}.
\end{proof}

The {\bf $K$-standardization map},
\[\Phi : \genomictableau{\nu / \lambda}{\mu} \to \inc{\nu / \lambda},\] 
is defined by filling the $k$th gene in the $<$-order with the entry $k$. 
Since any $T\in\genomictableau{\nu / \lambda}{\mu}$ is also a semistandard tableau
(by forgetting the gene structure) and since no two boxes of the same gene can be in the same row,
it follows that $\Phi(T)\in \inc{\nu / \lambda}$.
Note that when each gene is a single box, $\Phi$ is the usual standardization map.

\begin{example}
If $T$ is the genomic tableau \ 
$\begin{ytableau}
 *(lightgray)\blank & *(lightgray)\blank &*(SkyBlue) 1 &*(green) 2\\
*(red) 1 & *(SkyBlue) 1 &*(green) 2 \\
*(green) 2
\end{ytableau} \mbox{\ \ \ then $\Phi(T) = \ytableaushort{{*(lightgray)\blank} {*(lightgray)\blank} 23, 123,3}$.}
$ \qed
\end{example} 

A {\bf horizontal strip} is a skew shape with no two boxes in the same column.
Following \cite{Thomas.Yong:V}, a {\bf Pieri filling} is an increasing tableau of 
horizontal strip shape where, in addition, labels weakly 
increase from southwest to northeast. 

Let 
\begin{equation}\nonumber
{\mathcal P}_k(\mu):=\left\{ 1 + \sum_{i<k} \mu_i, 2 + \sum_{i<k} \mu_i, \dots, \sum_{j\leq k} \mu_j\right\}.
\end{equation}
That is, 
\[{\mathcal P}_1(\mu)=\{1,2,\ldots,\mu_1\}, \ {\mathcal P}_2(\mu)=\{\mu_1+1,\ldots,\mu_1+\mu_2\}, \text{\ etc.}\]
We say $S \in \inc{\nu / \lambda}$ is {\bf $\mu$-Pieri-filled} if for each $k \leq \ell(\mu)$, the entries of $S$ in
${\mathcal P}_k(\mu)$
form a Pieri filling of a horizontal strip. 

\begin{example}
The increasing tableau \!\!\!\!\!\!\! $\ytableaushort{\none \none 23, \none 14, 13}$ \, is not $(2, 2)$-Pieri-filled, as the entries 3 and 4 do not form a Pieri filling. However, it is $(2,1,1)$-Pieri-filled. \qed
\end{example}

Let 
\[\pieri{\nu / \lambda}{\mu}=\{\text{$S\in \inc{\nu/\lambda}$ that are $\mu$-Pieri-filled}\}.\]

\begin{theorem}\label{thm:Phi_Psi_bijection}
$\Phi:\genomictableau{\nu / \lambda}{\mu}\to \pieri{\nu / \lambda}{\mu}$ is a bijection.
\end{theorem}
\begin{proof}
We begin by defining the {\bf $K$-semistandardization map} 
\[\Psi:\pieri{\nu / \lambda}{\mu}\to \genomictableau{\nu / \lambda}{\mu}.\]
This extends the classical semistandardization map 
from standard Young tableaux to semistandard Young tableaux. Suppose $S\in {\tt PF}_{\mu}(\nu/\lambda)$. 
Construct a filling $T$ of $\nu/\lambda$ by placing into each box the unique positive integer $k$ such that $i \in \mathcal{P}_k(\mu)$, where $i$ is the entry of the corresponding box of $S$. 
Clearly, $T$ is a semistandard tableau.  

Declare boxes of $T$ to be in the same gene if and only if the corresponding boxes of $S$ contain the same value. 
Since $S$ is an increasing tableau, each gene of $T$ has at most one box in any row. 
Since the entries of $S$ in ${\mathcal P}_k(\mu)$
form a Pieri filling, given any two genes $\GG_1$, $\GG_2$ of family $k$ in $T$, 
every box $\GG_1$ appears west of every box of $\GG_2$ (or \emph{vice versa}).
Hence $T \in \genomictableau{\nu / \lambda}{\mu}$.

We now show that
$\Phi$ is well-defined, i.e.,
\[\im \Phi\subseteq {\tt PF}_{\mu}(\nu/\lambda).\]
Fix 
\[T \in \genomictableau{\nu / \lambda}{\mu}\] 
and set
\[S := \Phi(T) \in \inc{\nu / \lambda}.\] 
Let $\gamma\subseteq \nu/\lambda$ be the set of boxes that contain $k$ in $T$.
By (column) semistandardness of $T$, $\gamma$ is a horizontal strip.
Since $\Phi$ puts the labels of ${\mathcal P}_k(\mu)$ into $\gamma$ (in $S$) so as to increase southwest to northeast, 
the resulting filling is $\mu$-Pieri-filled.

By construction we have that 
\[\Phi\circ\Psi={\rm id}_{\pieri{\nu / \lambda}{\mu}} \text{\ and 
$\Psi\circ\Phi={\rm id}_{\genomictableau{\nu / \lambda}{\mu}}$.}\] 
It is straightforward from the definitions that
$\Psi$ are $\Phi$ are injective maps. Therefore, we conclude
that $\Phi$ and $\Psi$ are mutually inverse bijections.
\end{proof}

\section{Genomic words and Knuth equivalence}\label{sec:sequences}
A {\bf genomic word} is a word $s$ of colored positive integers such that all $i$'s of a fixed color are consecutive among the set of all $i$'s. A {\bf genotype} of $s$ is a subword that selects one letter of each color. Say $s$ is {\bf ballot} if every genotype of $s$ is ballot.

\begin{example}
${\color{green} 2} {\color{red} 1} {\color{green} 2} {\color{red} 1} {\color{SkyBlue} 1} {\color{green} 2}$ is a genomic word, whereas ${\color{green} 2} {\color{red} 1} {\color{green} 2} {\color{SkyBlue} 1} {\color{red} 1} {\color{green} 2}$ is not
because the subword of $1$'s is ${\color{red} 1}  {\color{SkyBlue} 1} {\color{red} 1}$ and the ${\color{red} 1}$'s are not 
consecutive.\qed
\end{example}

Let ${\tt genomicseq}(T)$ be the 
colored row reading word (taken in right to left and top to bottom order) of a genomic tableau $T$. 

\begin{lemma}
\label{lem:seqisword}
For a genomic tableau $T$, ${\tt genomicseq}(T)$ is a genomic word.
\end{lemma}
\begin{proof}
The follows from the semistandardness of $T$ together with the condition that the boxes of each gene of family $i$ are consecutive in the left-to-right order on $i$'s in $T$.
\end{proof}

We extend the $K$-standardization map $\Phi$ to genomic words by $\Phi(s) := \mathtt{seq}(\Phi(\widehat{T}(s)))$ where $\widehat{T}(s)$ is the antidiagonal of disconnected boxes filled from northeast to southwest by the given genomic word. 

\begin{lemma}\label{lem:word_stuff}
\gap
\begin{itemize}
\item[(I)]
Every genomic word $s$ is $\mathtt{genomicseq}(T)$ for some genomic tableau $T$.
\item[(II)] $T$ is ballot if and only if $\mathtt{genomicseq}(T)$ is ballot.
\item[(III)] If $\mathtt{genomicseq}(T) = s$, then $\Phi(s) = \mathtt{seq}(\Phi(T))$. 
\end{itemize}
\end{lemma}
\begin{proof}
For (I), in particular, one can take $T = \widehat{T}(s)$. 
(II) is clear. (III) is straightforward.
\end{proof}

\begin{example}
If $T$ is the genomic tableau 
$\begin{ytableau}
 *(lightgray)\blank & *(lightgray)\blank &*(red) 1 &*(green) 2\\
*(cyan) 1 & *(red) 1 &*(green) 2 \\
*(green) 2
\end{ytableau}$, then $\mathtt{genomicseq}(T) = {\color{green} 2} {\color{red} 1} {\color{green} 2} {\color{red} 1} {\color{SkyBlue} 1} {\color{green} 2}$.
By selecting one {\color{green} green} letter, one {\color{red} red} letter, and one {\color{SkyBlue} blue} letter from ${\color{green}2}{\color{red}1}{\color{green}2}{\color{red}1}{\color{SkyBlue} 1}{\color{green}2}$ we arrive at three possible genotypes
of $\mathtt{genomicseq}(T)$: $211$, $121$ and $112$. Thus $\mathtt{genomicseq}(T)$ is not ballot. \qed
\end{example}

{\bf Genomic Knuth equivalence} is the equivalence relation $\equiv_{G}$ on genomic words obtained as the transitive closure of
\begin{align}
\mathbf{u}{\color{red}{ii}}\mathbf{v} &\equiv_{G} \mathbf{u}{\color{red}{i}}\mathbf{v},\tag{G.1} \\
\mathbf{u}{\color{red}{i}}{\color{SkyBlue}{j}}{\color{red}{i}}\mathbf{v} &\equiv_{G} \mathbf{u}{\color{SkyBlue}{j}}{\color{red}{i}}{\color{cyan}{j}}\mathbf{v},\tag{G.2} \\ 
\mathbf{u}{\color{SkyBlue}{j}}{\color{red}{i}}{\color{green}{k}}\mathbf{v} &\equiv_{G} \mathbf{u}{\color{SkyBlue}{j}}{\color{green}{k}}{\color{red}{i}}\mathbf{v},\tag{G.3} \\ 
\mathbf{u}{\color{red}{p}}{\color{green}{q}}{\color{SkyBlue}{j}}\mathbf{v} &\equiv_{G} \mathbf{u}{\color{green}{q}}{\color{red}{p}}{\color{SkyBlue}{j}}\mathbf{v},\tag{G.4}
\end{align}
where $i \leq j < k$, $p < j \leq q$, and {\color{red}{red}}, {\color{SkyBlue}{blue}}, {\color{green}{green}} are distinct colors. This equivalence relation is a genomic version of the \emph{$K$-Knuth equivalence} introduced by A.~Buch--M.~Samuel \cite[$\mathsection 5$]{Buch.Samuel}. It furthermore generalizes Knuth equivalence \cite{Knuth} in the sense that it agrees with this older notion on words where each letter is of a distinct color, obviating (G.1) and (G.2).

\begin{theorem}\label{thm:knuth_preserves_ballot}
If ${\bf x} \equiv_{G} {\bf y}$, then ${\bf x}$ is ballot if and only if ${\bf y}$ is ballot.
\end{theorem}
\begin{proof}
Let ${\bf x}$ be a genomic word. It suffices to show that (G.1)--(G.4)  do not change the ballotness of ${\bf x}$.

(G.1) and (G.2) preserve the set of genotypes and therefore ballotness.

(G.3) clearly preserves ballotness unless $k = i + 1$.
In this case, since $i\leq j<k$, this means $i=j$.
Suppose therefore 
\[{\bf x} = \mathbf{u}{\color{SkyBlue}{j}}{\color{red}{j}}{\color{green}{k}}\mathbf{v}, \mbox{\ and that ${\bf y} = \mathbf{u}{\color{SkyBlue}{j}}{\color{green}{k}}{\color{red}{j}}\mathbf{v}$.}\] 

Clearly if ${\bf x}$ is not ballot, then ${\bf y}$ is not ballot. Conversely,
assume ${\bf x}$ is ballot. 
It is enough to show that $\mathbf{u}{\color{SkyBlue}{j}}{\color{green}{k}}$ is ballot. Since ${\bf x}$ is ballot, the initial segment $\mathbf{u}{\color{SkyBlue}{j}}$ is ballot. Now, deleting the last letter of a ballot word leaves a ballot word. Since the last letter in the case at hand is
${\color{SkyBlue}{j}}$ it follows that 
the subsequence of $\mathbf{u}{\color{SkyBlue}{j}}$ formed by deleting every ${\color{SkyBlue}{j}}$ is ballot. In particular, every genotype of $\mathbf{u}{\color{SkyBlue}{j}}$ has strictly more $j$'s than $k$'s. Thus $\mathbf{u}{\color{SkyBlue}{j}}{\color{green}{k}}$ (and hence ${\bf y}$) is ballot.

(G.4) is only a concern if $q=p+1$. In this case, since $p<j\leq q$, we must also have $j=q$. Thus, 
\[{\bf x} = \mathbf{u}{\color{red}{p}}{\color{green}{q}}{\color{cyan}{q}}\mathbf{v} \text{ and ${\bf y} = \mathbf{u}{\color{green}{q}}{\color{red}{p}}{\color{SkyBlue}{q}}\mathbf{v}$.}\] 
If ${\bf y}$ is ballot, then ${\bf x}$ is ballot. 
Conversely, assume ${\bf x}$ is ballot. It suffices to show that $\mathbf{u}{\color{green}{q}}{\color{red}{p}}$ is ballot. Since it is an initial segment of ${\bf x}$, $\mathbf{u}{\color{red}{p}}{\color{green}{q}}{\color{SkyBlue}{q}}$ is ballot. Given any two genes of family $q$ in any genomic word, one appears entirely right of the other. Thus ${\color{SkyBlue}{q}}$ does not appear in $\mathbf{u}$, and hence every genotype of $\mathbf{u}{\color{red}{p}}{\color{green}{q}}$ has strictly more $p$'s than $q$'s. Thus $\mathbf{u}{\color{green}{q}}{\color{red}{p}}$ is ballot.
\end{proof}

\section{Genomic jeu de taquin}\label{sec:genomic_jdt}

If $T \in \genomictableau{\nu / \lambda}{\mu}$, an {\bf inner corner} of $T$ is a maximally southeast box of $\lambda$. Let $I$ be any set of inner corners of $T$. We obtain a genomic tableau $\jdt{T}{I}$ as follows: Place a $\bullet$ in each box of $I$; let $T^{\bullet}$ denote the result. Two boxes of a tableau are {\bf neighbors} if they share a 
horizontal or vertical edge. For each gene $\GG$, define the
operator $\mathtt{switch}_{\GG}^{\bullet}$ as follows. 
Every box of $\GG$ with a 
neighbor containing a $\bullet$ 
becomes a box containing a $\bullet$, while simultaneously every box with a $\bullet$ and a $\GG$ 
neighbor becomes a box of $\GG$. The remaining boxes are unchanged by $\mathtt{switch}_{\GG}^{\bullet}$. 

Index the genes of $T$ as \[\GG_1 < \GG_2 < \dots < \GG_{|\mu|}\] 
according to the total order on genes from Lemma~\ref{lem:gene_total_order}.  Then 
\[\jdt{T}{I}:=\mathtt{switch}_{\GG_{|\mu|}}^{\bullet} \circ \dots \circ \mathtt{switch}_{\GG_2}^{\bullet} \circ \mathtt{switch}_{\GG_1}^{\bullet}(T^{\bullet})\] 
with the $\bullet$'s deleted. (This algorithm reduces to M.-P.~Sch\"utzenberger's \emph{jeu de taquin} for semistandard tableaux in the case each gene contains only a single box.)

\begin{example} Suppose $T^\bullet$ is the genomic tableau
$\begin{ytableau}
*(lightgray)\blank & *(lightgray)\blank & *(lightgray)\blank & \bullet \\
 *(lightgray)\blank & \bullet &*(red) 1 &*(green) 2\\
*(SkyBlue) 1 & *(red) 1 &*(green) 2 \\
*(green) 2
\end{ytableau}$. Then
\[\mathtt{switch}^\bullet_{{\color{SkyBlue}1}}(T^\bullet) = 
\begin{ytableau}
*(lightgray)\blank & *(lightgray)\blank & *(lightgray)\blank & \bullet \\
 *(lightgray)\blank & \bullet & *(red) 1 & *(green) 2\\
*(SkyBlue) 1 & *(red) 1 & *(green) 2 \\
*(green) 2
\end{ytableau}, \ \ \ 
\mathtt{switch}_{\color{red}1}^{\bullet} \circ \mathtt{switch}_{\color{SkyBlue}1}^{\bullet}(T^\bullet) = \begin{ytableau}
*(lightgray)\blank & *(lightgray)\blank & *(lightgray)\blank & \bullet \\
 *(lightgray)\blank & *(red) 1 & \bullet &*(green) 2\\
*(SkyBlue) 1 & \bullet &*(green) 2 \\
*(green) 2
\end{ytableau},\ \ \text{and }\]
\[\mathtt{switch}_{\color{green}2}^{\bullet} \circ \mathtt{switch}_{\color{red}1}^{\bullet} \circ \mathtt{switch}_{\color{SkyBlue}1}^{\bullet}(T^\bullet) = \begin{ytableau}
*(lightgray)\blank & *(lightgray)\blank & *(lightgray)\blank & *(green) 2 \\
 *(lightgray)\blank & *(red) 1 & *(green) 2 &\bullet\\
*(SkyBlue) 1 & *(green) 2 &\bullet \\
*(green) 2
\end{ytableau}.\ \
\text{So } \jdt{T}{I} = \begin{ytableau}
*(lightgray)\blank & *(lightgray)\blank & *(lightgray)\blank & *(green) 2 \\
 *(lightgray)\blank & *(red) 1 & *(green) 2\\
*(SkyBlue) 1 & *(green) 2 \\
*(green) 2
\end{ytableau}.\] \qed
\end{example}

Define {\bf jeu de taquin equivalence} $\sim_G$ on genomic tableaux as the symmetric, transitive closure of the relation $T \sim_G \jdt{T}{I}$. We now state the genomic analogue of \cite[Theorem~6.2]{Buch.Samuel}
(restated as Theorem~\ref{thm:BuchSamuel} below):

\begin{theorem}\label{thm:genomic tableau_knuth_equivalence_determines_plactic_class}
Let $T, U$ be genomic tableaux. Then $T \sim_G U$ if and only if $\mathtt{genomicseq}(T) \equiv_{G} \mathtt{genomicseq}(U)$.
\end{theorem}
\begin{proof}
We assume the terminology and results on the ${\tt Kjdt}$ (\emph{jeu de taquin}) for 
increasing tableaux of
\cite[$\mathsection$1]{Thomas.Yong:V} (see specifically pages 123--124). {\bf $K$-Knuth equivalence} \cite[$\mathsection$5]{Buch.Samuel} is the 
the symmetric, transitive closure of the following {\bf $K$-Knuth relations} 
(our conventions are reversed from those of \cite{Buch.Samuel}; this has no effect on the applicability of their results):
For words ${\bf u}, {\bf v}$ and integers $0<i<j<k$,
\begin{align}
\mathbf{u}{{ii}}\mathbf{v} &\equiv_{K} \mathbf{u}{{i}}\mathbf{v},\tag{K.1} \\
\mathbf{u}{{i}}{{j}}{{i}}\mathbf{v} &\equiv_{K} \mathbf{u}{{j}}{{i}}{{j}}\mathbf{v},\tag{K.2} \\ 
\mathbf{u}{{j}}{{i}}{{k}}\mathbf{v} &\equiv_{K} \mathbf{u}{{j}}{{k}}{{i}}\mathbf{v},\tag{K.3} \\ 
\mathbf{u}{{i}}{{k}}{{j}}\mathbf{v} &\equiv_{K} \mathbf{u}{{k}}{{i}}{{j}}\mathbf{v}.\tag{K.4}
\end{align}

Define {\bf ${\tt Kjdt}$-equivalence} ($\sim_K$) on increasing tableaux as the symmetric, transitive closure of the relation $T \sim_K {\tt Kjdt}_I(T)$.
The key relationship between these two equivalence relations is:

\begin{theorem}\cite[Theorem~6.2]{Buch.Samuel}\label{thm:BuchSamuel}
 $T \sim_K U$ if and only if ${\tt seq}(T) \equiv_K {\tt seq}(U)$. \qed
\end{theorem}

Let ${\tt Gen}(\nu/\lambda)$ be the set of all genomic tableaux 
of shape $\nu /\lambda$.

\begin{lemma}\label{lem:genomic tableau_jdt}
For $T\in {\tt Gen}(\nu/\lambda)$ and set of inner corners $I$, $\Phi(\jdt{T}{I}) = {\tt Kjdt}_I({\Phi(T)})$.
\end{lemma}
\begin{proof}
From the definitions, this is an easy induction on the number of genes of $T$.
\end{proof}

\begin{lemma}
\label{lem:wordequivsame}
For any genomic words ${\bf u}$ and ${\bf v}$, we have
${\bf u} \equiv_{G} {\bf v}$ if and only if
$\Phi({\bf u}) \equiv_{K} \Phi({\bf v})$.
\end{lemma}
\begin{proof}
Immediate from the definitions of $\equiv_{K}$ and $\equiv_{G}$.
\end{proof}

By Lemma~\ref{lem:genomic tableau_jdt}, $T \sim_G U$ if and only if $\Phi(T) \sim_K \Phi(U)$. By Theorem~\ref{thm:BuchSamuel}, the latter relation is equivalent to 
\[\mathtt{seq}(\Phi(T)) \equiv_K \mathtt{seq}(\Phi(U)).\] 
By Lemma~\ref{lem:wordequivsame} and Lemma~\ref{lem:word_stuff}(III), we see that
\[\mathtt{seq}(\Phi(T)) \equiv_K \mathtt{seq}(\Phi(U))\] 
is equivalent to 
\[\mathtt{genomicseq}(T) \equiv_G \mathtt{genomicseq}(U).\]
\end{proof}

\begin{corollary}\label{cor:genomic tableau_jdt_preserves_ballotness}
If $T\sim_G U$, then $T$ is ballot if and only if $U$ is ballot.
\end{corollary}
\begin{proof}
By Theorem~\ref{thm:genomic tableau_knuth_equivalence_determines_plactic_class} and Theorem~\ref{thm:knuth_preserves_ballot}.
\end{proof}

Let $T_{\mu}$ be the {\bf highest weight tableau} of shape $\mu$,
i.e., the semistandard tableau whose $i$-th row uses only the label $i$. Note $T_{\mu}$ may be
also regarded as a genomic tableau in a unique manner. Let $S_{\mu}:=\Phi(T_{\mu})$
be the {\bf row superstandard tableau} of shape $\mu$ (this is the tableau whose first row has entries $1,2,3,\ldots,\mu_1$, and whose
second row has entries $\mu_1+1,\mu_2+2,\ldots,\mu_1+\mu_2$ etc.). 

\begin{corollary}[of Lemma~\ref{lem:genomic tableau_jdt}]
\label{cor:gettouse}
For $T \in {\tt Gen}(\nu / \lambda)$, $T \sim_G T_\mu$ if and only if
$\Phi(T) \sim_K S_\mu$. 
\end{corollary}
\begin{proof}
This is immediate from Lemma~\ref{lem:genomic tableau_jdt}
because $S_\mu = \Phi(T_\mu)$. 
\end{proof}

\section{Three proofs of the Genomic Littlewood-Richardson rule (Theorem \ref{thm:A_lr_rule})}

\subsection{Proof 1: Bijection with increasing tableaux}\label{sec:first_proof}
Our first proof uses the results of Sections~\ref{sec:K-maps}--\ref{sec:genomic_jdt} to prove
Theorem \ref{thm:A_lr_rule}. Let 
\[\ballot{\nu / \lambda}{\mu}:=\{\text{$T \in \genomictableau{\nu / \lambda}{\mu}:$ $T$ is ballot}\}.\]
Also, let 
\[\incr{\nu / \lambda}{\mu}:=
\{T\in \inc{\nu/\lambda}: {\tt KRect}(T)=S_\mu\}.\]

\begin{lemma}\label{lem:thefirst}
Let $T \in \genomictableau{\nu / \lambda}{\mu}$.
Then $T\in \ballot{\nu/\lambda}{\mu}$ if and only if 
$\Phi(T)\in \incr{\nu/\lambda}{\mu}$.
\end{lemma}
\begin{proof}
Suppose $T$ is ballot. By iterating application of ${\tt jdt}_I$ (under arbitrary choices of nonempty sets $I$ of inner corners) starting with $T$, we have that
$T\sim_G R$ for some straight-shaped tableau $R$ ({\it a priori}, $R$ might depend on the choices of $I$). By Corollary~\ref{cor:genomic tableau_jdt_preserves_ballotness}, $R$ is ballot. 
Since genomic {\it jeu de taquin} preserves tableau content, $R=T_{\mu}$.
Hence, by Lemma~\ref{lem:genomic tableau_jdt}, $\Phi(T)$ rectifies to $S_\mu$.

Conversely, suppose $\Phi(T)$ rectifies to $S_\mu$. Then by Lemma~\ref{lem:genomic tableau_jdt}, $T$ rectifies to $T_{\mu}$. But $T_{\mu}$ is a ballot genomic tableau. Hence by Corollary~\ref{cor:genomic tableau_jdt_preserves_ballotness}, $T$ is also ballot.
\end{proof}

\begin{lemma}\label{lem:superstandard_rectifiers_are_Pieri_filled}
$\incr{\nu / \lambda}{\mu}\subseteq \pieri{\nu/\lambda}{\mu}$.
\end{lemma}
\begin{proof}
This is part of \cite[Proof of Theorem 1.2]{Thomas.Yong:V}. 
\end{proof}

In view of Lemmas~\ref{lem:thefirst} and~\ref{lem:superstandard_rectifiers_are_Pieri_filled}, we may define 
\[\phi : \ballot{\nu / \lambda}{\mu} \to \incr{\nu / \lambda}{\mu}\] 
as the restriction 
\[\Phi|_{\ballot{\nu / \lambda}{\mu}}\] and define 
\[\psi : \incr{\nu / \lambda}{\mu} \to \ballot{\nu / \lambda}{\mu}\]
as the restriction 
\[\Psi|_{\incr{\nu / \lambda}{\mu}}.\] 
Now $\Phi$ and $\Psi$ are
mutually inverse bijections (cf.~Theorem~\ref{thm:Phi_Psi_bijection}). Thus
$\phi$ and $\psi$ are mutually inverse bijections between $\incr{\nu / \lambda}{\mu}$ and $\ballot{\nu / \lambda}{\mu}$. Hence the theorem follows from
the ${\tt Kjdt}$ rule of \cite{Thomas.Yong:V} for $a_{\lambda,\mu}^{\nu}$.\qed

\subsection{Proof 2: Bijection with set-valued tableaux}\label{sec:bijection_with_Buch}

In our next proof, we relate genomic tableaux to
the original rule for $a_{\lambda,\mu}^{\nu}$ found by A.~Buch \cite[Theorem~5.4]{Buch:KLR}.

We first recall some definitions from \cite{Buch:KLR}. 
A {\bf set-valued tableau} $T$ of (skew) shape
$\nu / \lambda$ is a filling of the boxes of $\nu / \lambda$ with non-empty finite subsets of ${\mathbb N}$ with the property that any tableau obtained by choosing exactly one label from each box is a (classical) semistandard tableau. The {\bf column reading word} of $T$, denoted ${\tt colword}(T)$ is obtained by reading the entries
of $T$ from bottom to top along columns and from left to right.
The entries in a non-singleton box are read in increasing order.
Such a word $(w_1,w_2,\ldots,w_N)$ is a {\bf reverse lattice word} if the content of $(w_L,w_{L+1},\ldots,w_N)$ is a partition for
every $1\leq L\leq N$, that is to say if its reverse is ballot. Finally, the shape $\mu\star \lambda$ is the skew shape
obtained by placing $\mu$ and $\lambda$ in southwest to northeast orientation with $\mu$'s northeast corner incident
to $\lambda$'s southwest corner. In other words 
\[\mu\star\lambda=(\mu_1+\lambda_1,\ldots,\mu_1+\lambda_{\ell(\lambda)},\mu_1,\mu_2,\ldots,\mu_{\ell(\mu)}) / (\mu_1^{\ell(\lambda)}).\]

\begin{example}
\ytableausetup{boxsize = 0.4em,aligntableaux=center}
If $\lambda = \ydiagram[*(red)]{2,2}$ and $\mu = \ydiagram[*(SkyBlue)]{3,1}$, then $\mu \star \lambda = \ydiagram{3+2,3+2,3,1} *[*(SkyBlue)]{0,0,3,1} *[*(red)]{3+2,3+2}$.
\qed
\ytableausetup{boxsize = 1.1em}
\end{example}

\begin{theorem}[{A.~Buch \cite[Theorem~5.4]{Buch:KLR}}]
\label{thm:buchsrule}
$(-1)^{|\lambda|+|\mu|-|\nu|}a_{\lambda,\mu}^{\nu}$ equals the number of set-valued tableaux $T$ of shape $\mu\star\lambda$ and content $\nu$
such that ${\tt colword}(T)$ is reverse lattice.
\end{theorem}

Let ${\tt Buch}_\nu(\mu \star \lambda)$ be the set of tableaux from Theorem~\ref{thm:buchsrule}. We define a map
\[  \Xi :{\tt Buch}_\nu(\mu \star \lambda) \to \ballot{\nu / \lambda}{\mu}\]
as follows. Let $T\in{\tt Buch}_\nu(\mu \star \lambda)$. Start with an empty shape $\lambda$. Read the columns of the $\mu$ portion of $T$ from top to bottom and right to left. Suppose a set 
\[S=\{s_1<\ldots<s_t\}\]
gives the entries of a box in row $i$ of the $\mu$ shape. Then place a new gene of family $i$ in the rows $s_1,\ldots,s_t$ (as far left
as possible in each case). Then $\Xi$ clearly has a (putative) inverse 
\[\Theta:\ballot{\nu / \lambda}{\mu}\to {\tt Buch}_\nu(\mu \star \lambda)\]
that records in row $i$ and column $1$ of the $\mu$ shape the rows that the leftmost gene of family $i$ sits in.
Similarly, in row $i$ and column $2$ we record the rows that the second leftmost gene of family $i$ sits in, etc.

\begin{theorem}
\label{theorem:BuchBallotbij}
$\Xi:{\tt Buch}_\nu(\mu \star \lambda)\to \ballot{\nu / \lambda}{\mu}$ and 
$\Theta:\ballot{\nu / \lambda}{\mu}\to {\tt Buch}_\nu(\mu \star \lambda)$ are well-defined and mutually inverse
bijections.
\end{theorem}

\begin{example}\label{ex:Buch}
\ytableausetup{boxsize=1.5em}
Let $\lambda = (2,1), \mu = (1,1)$ and $\nu = (3,2,1)$. Then ${\tt Buch}_\nu(\mu \star \lambda)$ consists of the two tableaux
\[ B_1 = \ytableaushort{\none  1 1, \none 2, {1,2}, 3 }  \; \text{and } B_2 = \ytableaushort{\none  1 1, \none 2, 1, {2,3} }. \]
We have 
\ytableausetup{boxsize=1.1em}
\[ \Xi(B_1) = \ytableaushort{\blank \blank {*(SkyBlue) 1}, \blank {*(SkyBlue) 1}, {*(red) 2}} \; \text{and } \Xi(B_2) = \ytableaushort{\blank \blank {*(SkyBlue) 1}, \blank {*(red) 2}, {*(red) 2}} .\]
The reader can check that these are the unique two elements of $\ballot{\nu / \lambda}{\mu}$.
\qed
\ytableausetup{boxsize=1.1em}
\end{example}

\begin{proof}[Proof of Theorem~\ref{theorem:BuchBallotbij}]
Let $T\in {\tt Buch}_\nu(\mu \star \lambda)$ and set $U:=\Xi(T)$.

\noindent
({\sf ${\Xi}$ is well-defined}): By definition, the number of genes of family $i$ is $\mu_i$. Hence the content of $U$ is $\mu$, as required. Next, observe that
since in each row of $T$ the entries increase weakly from left to right, no two genes of the same family interweave.
Also note that no two labels of the same gene are in the same column since otherwise we would obtain that ${\tt colword}(T)$ is not reverse lattice, since labels in the same box are
read in increasing order, a contradiction.

The hypothesis that ${\tt colword}(T)$ is reverse lattice precisely guarantees that when adding the boxes in the rows of $S$ 
one takes a Young diagram to a larger Young diagram. Thus $U$ is a tableau of (skew) Young diagram shape. Note that since $S$ is a set,
no row of $U$ contains two boxes of the same gene. 

We next verify the semistandardness conditions. 
Suppose $U$ violates the horizontal semistandardness requirement.
That is, there is a box ${\sf x}$ directly left and adjacent to a box ${\sf y}$ in $U$ such that ${\tt lab}_U({\sf x})>{\tt lab}_U({\sf y})$.
Let ${\sf x}'$ and ${\sf y}'$ be the boxes in $T$ that added ${\sf x}$ and ${\sf y}$ during the execution of $\Xi$. Since 
${\tt lab}_U({\sf x})>{\tt lab}_U({\sf y})$, by $\Xi$'s definition, the row of ${\sf x}'$ is strictly south of the row of ${\sf y}'$.
Moreover, since ${\sf x}$ is left
of ${\sf y}$ we know that ${\sf x}'$ is read before ${\sf y}'$ in ${\tt colword}(T)$. Therefore, ${\sf y}'$ is strictly north and strictly
west of ${\sf x}'$. However, since $T$ is a (set-valued) semistandard tableau, the labels of ${\sf y}'$ in $T$ are all strictly smaller than 
those of ${\sf x}'$. This implies that ${\sf y}$ is in a row strictly north of that of ${\sf x}$, a contradiction.
The argument that $U$ satisfies the vertical semistandardness requirement is similar.

It remains to check that $U$ is ballot. To do this, make an arbitrary but fixed choice
of genotype $G_U$ of $U$.
The labels of family $i$ and $i+1$ may be blamed
on labels in rows $i$ and $i+1$ of $T$. Suppose the sets of labels in those rows are  
\[Q_1,Q_2,\ldots,Q_t, Q_{t+1},\ldots,Q_{t+s} \text{\ \ (row $i$)
and  \ $R_1,R_2,\ldots,R_t$ \ (row $i+1$)}\] 
where $s\geq 0$.
Since we know $U$ is semistandard, the labels associated to rows $i$ and $i+1$ separately form a Pieri strip.  Here $Q_1$ is associated to the rightmost gene of family $i$ (in $U$) and $Q_{t+s}$ is associated to the leftmost gene of family $i$ (in $U$). Similarly, $R_1$ is
associated to the rightmost gene of family $i+1$ (in $U$) 
and $R_{t}$ is associated to the leftmost gene of family $i+1$ (in $U$). By the vertical
semistandardness of $T$, we have 
\[\max Q_i< \min R_i \text{\ for $1\leq i\leq t$.}\]
This clearly implies  that the $m$th rightmost label of family $i+1$ in $G_U$ is strictly south and weakly west of the $m$th rightmost label of family $i$ in $G_U$, for $1\leq m\leq t$. Since this is true for each $i$, $G_U$ is ballot.

\noindent
({\sf ${\Xi}$ is injective}): Clear.

\noindent
({\sf $\Theta$ is well-defined}:) This is proved with the same
arguments (said in reverse) as those given in the well-definedness of $\Xi$.

\noindent
({\sf $\Theta$ is injective}): Clear.

The theorem follows since $\Xi$ and $\Theta$ are mutually inverse injections.
\end{proof}

Composing Theorem~\ref{theorem:BuchBallotbij}  with the bijection of Section~\ref{sec:first_proof} permits one to biject the above rule of A.~Buch
with the $K$-theoretic \emph{jeu de taquin} rule of \cite{Thomas.Yong:V}.

\subsection{Proof 3: Bijection with puzzles}
A third proof of Theorem~\ref{thm:A_lr_rule} considers the bijection given in \cite{PY:puzzles} between more general genomic tableaux and the Knutson-Vakil puzzles of \cite[$\mathsection 5$]{Coskun.Vakil}. It is not hard to see that 
this restricts to a bijection between the tableaux of Theorem~\ref{thm:A_lr_rule} and the ordinary $K$-theory puzzles of A.~Buch \cite[$\mathsection 3.3$]{Vakil}. Since the latter are known to calculate $a_{\lambda, \mu}^\nu$, Theorem~\ref{thm:A_lr_rule} follows. 

Consider the $n$-length equilateral triangle oriented as $\Delta$. A {\bf puzzle}
is a filling of $\Delta$ with the following {\bf puzzle pieces}:
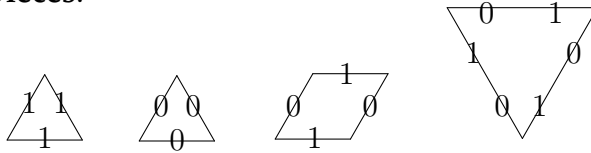
\[\begin{picture}(250,45)
\put(0,0){
\begin{tikzpicture}[line cap=round,line join=round,>=triangle 45,x=1.0cm,y=1.0cm]
\clip(0.7,-0.38) rectangle (2.22,1.08);
\draw (1,0)-- (2,0);
\draw (2,0)-- (1.5,0.87);
\draw (1.5,0.87)-- (1,0);
\end{tikzpicture}}
\put(17,21){$1$}
\put(29,21){$1$}
\put(23,8){$1$}

\put(50,4){
\begin{tikzpicture}[line cap=round,line join=round,>=triangle 45,x=1.0cm,y=1.0cm]
\clip(0.7,-0.23) rectangle (2.26,0.97);
\draw (1,0)-- (2,0);
\draw (2,0)-- (1.5,0.87);
\draw (1.5,0.87)-- (1,0);
\end{tikzpicture}}
\put(67,20){$0$}
\put(79,20){$0$}
\put(73,7){$0$}

\put(100,1){
\begin{tikzpicture}[line cap=round,line join=round,>=triangle 45,x=1.0cm,y=1.0cm]
\clip(0.65,-0.36) rectangle (2.86,1.08);
\draw  (1,0)-- (1.5,0.87);
\draw (1.5,0.87)-- (2.5,0.87);
\draw  (2.5,0.87)-- (2,0);
\draw (2,0)-- (1,0);
\end{tikzpicture}
}
\put(117,20){$0$}
\put(146,20){$0$}
\put(125,8){$1$}
\put(137,32){$1$}

\put(140,4){
\begin{tikzpicture}[line cap=round,line join=round,>=triangle 45,x=1.0cm,y=1.0cm]
\clip(-0.23,-2) rectangle (3,1);
\draw (1,0)-- (3,0);
\draw (3,0)-- (2,-1.74);
\draw (2,0-1.74)-- (1,0);
\end{tikzpicture}}
\put(190,55){$0$}
\put(216,55){$1$}
\put(196,20){$0$}
\put(210,20){$1$}
\put(185,40){$1$}
\put(223,40){$0$}
\end{picture}\]
Henceforth, we color code these pieces as black, white, gray, and blue respectively, dropping the numerical labels.
A {\bf filling} requires that 
the common edges of adjacent puzzle pieces share the same label. 
The first three may be rotated but the fourth
({\bf $K$-piece}) may not.  A {\bf $K$-puzzle} is a puzzle filling of $\Delta$. 

Convert partitions inside a $k \times (n-k)$ rectangle to binary sequences in the following way. Starting at the upper right corner of the $k \times (n-k)$ rectangle, construct a lattice path from the binary sequence by reading each segment $(-1,0)$ as $0$ and each segment $(0,-1)$ as $1$. For example, we convert the partition $(3,2)$ to the binary string $010100$ as follows:
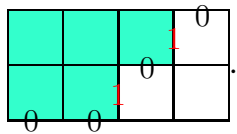
\[
\begin{picture}(120,37)
\ytableausetup{boxsize=1.7em}
\put(0,18){$\ytableaushort{ {*(ttffcc)\blank} {*(ttffcc)\blank} {*(ttffcc)\blank} {*(white)\blank}, 
{*(ttffcc)\blank} {*(ttffcc)\blank} {*(white)\blank} {*(white)\blank} }$.}
\put(39,6){$\textcolor{red} 1$}
\put(60,27){$\textcolor{red} 1$}
\put(6,-4){$0$}
\put(30,-4){$0$}
\put(50,16){$0$}
\put(71,36){$0$}
\end{picture} \]
Let $\Delta_{\lambda,\mu,\nu}$ be $\Delta$ with the boundary given by binary sequences
\begin{itemize}
\item $\lambda$ as read $\nearrow$ along the left side;
\item $\mu$ as read $\searrow$ along the right side; and
\item $\nu$ as read $\rightarrow$ along the bottom side.
\end{itemize}

\begin{theorem}[{A.~Buch \cite[$\mathsection 3.3$]{Vakil}}]
$(-1)^{|\nu| - |\lambda| - |\mu|} a_{\lambda,\mu}^{\nu}= \# \{ \text{$K$-puzzles of $\Delta_{\lambda,\mu,\nu}$} \}$. \qed
\end{theorem}

\begin{example}
Continuing Example~\ref{ex:Buch} and assuming the Grassmannian in question is ${\rm Gr}_3({\mathbb C}^6)$ , the bijection of \cite{PY:puzzles} matches $\Xi(B_1)$ and $\Xi(B_2)$ to the puzzles
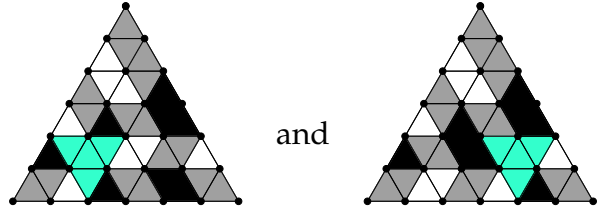
\[\Scale[0.5]{
\begin{tikzpicture}[line cap=round,line join=round,>=triangle 45,x=1.0cm,y=1.0cm]
\clip(-1.,-1.) rectangle (7.,6.);
\fill[color=ffffff,fill=ffffff,fill opacity=0.1] (0.,0.) -- (6.,0.) -- (3.,5.196152422706633) -- cycle;
\fill[color=aqaqaq,fill=aqaqaq,fill opacity=1.0] (0.,0.) -- (1.,0.) -- (1.5,0.8660254037844384) -- (0.5,0.8660254037844386) -- cycle;
\fill[color=aqaqaq,fill=aqaqaq,fill opacity=1.0] (6.,0.) -- (5.5,0.8660254037844439) -- (4.5,0.8660254037844394) -- (5.,0.) -- cycle;
\fill[fill=black,fill opacity=1.0] (0.5,0.8660254037844386) -- (1.,1.7320508075688772) -- (1.5,0.8660254037844384) -- cycle;
\fill[color=aqaqaq,fill=aqaqaq,fill opacity=1.0] (3.,5.196152422706633) -- (2.5,4.330127018922193) -- (3.,3.4641016151377544) -- (3.5,4.330127018922193) -- cycle;
\fill[color=aqaqaq,fill=aqaqaq,fill opacity=1.0] (3.5,4.330127018922193) -- (4.,3.4641016151377526) -- (3.5,2.598076211353316) -- (3.,3.4641016151377544) -- cycle;
\fill[fill=black,fill opacity=1.0] (4.,3.4641016151377526) -- (3.5,2.598076211353316) -- (4.5,2.598076211353314) -- cycle;
\fill[color=aqaqaq,fill=aqaqaq,fill opacity=1.0] (2.,3.4641016151377544) -- (1.5,2.598076211353316) -- (2.,1.732050807568877) -- (2.5,2.598076211353316) -- cycle;
\fill[fill=black,fill opacity=1.0] (4.5,0.8660254037844394) -- (4.,0.) -- (5.,0.) -- cycle;
\fill[color=ttffcc,fill=ttffcc,fill opacity=1.0] (1.,1.7320508075688772) -- (3.,1.7320508075688774) -- (2.,0.) -- cycle;
\fill[fill=black,fill opacity=1.0] (2.5,0.8660254037844379) -- (2.,0.) -- (3.,0.) -- cycle;
\fill[fill=black,fill opacity=1.0] (2.5,2.598076211353316) -- (3.,1.7320508075688774) -- (2.,1.732050807568877) -- cycle;
\fill[color=aqaqaq,fill=aqaqaq,fill opacity=1.0] (2.5,2.598076211353316) -- (3.5,2.598076211353316) -- (4.,1.732050807568878) -- (3.,1.7320508075688774) -- cycle;
\fill[fill=black,fill opacity=1.0] (3.5,2.598076211353316) -- (4.,1.732050807568878) -- (4.5,2.598076211353314) -- cycle;
\fill[fill=black,fill opacity=1.0] (4.,1.732050807568878) -- (5.,1.732050807568878) -- (4.5,2.598076211353314) -- cycle;
\fill[color=aqaqaq,fill=aqaqaq,fill opacity=1.0] (4.,1.732050807568878) -- (3.5,0.8660254037844379) -- (4.5,0.8660254037844394) -- (5.,1.732050807568878) -- cycle;
\fill[fill=black,fill opacity=1.0] (3.5,0.8660254037844379) -- (4.,0.) -- (4.5,0.8660254037844394) -- cycle;
\fill[color=aqaqaq,fill=aqaqaq,fill opacity=1.0] (3.5,0.8660254037844379) -- (2.5,0.8660254037844379) -- (3.,0.) -- (4.,0.) -- cycle;
\draw (0.,0.)-- (6.,0.);
\draw [color=ffffff] (0.,0.)-- (6.,0.);
\draw [color=ffffff] (6.,0.)-- (3.,5.196152422706633);
\draw [color=ffffff] (3.,5.196152422706633)-- (0.,0.);
\draw (0.5,0.8660254037844386)-- (1.,0.);
\draw (1.,1.7320508075688772)-- (2.,0.);
\draw (1.5,2.598076211353316)-- (3.,0.);
\draw (2.,3.4641016151377544)-- (4.,0.);
\draw (5.,0.)-- (2.5,4.330127018922193);
\draw (3.5,4.330127018922193)-- (1.,0.);
\draw (2.,0.)-- (4.,3.4641016151377526);
\draw (4.5,2.598076211353314)-- (3.,0.);
\draw (4.,0.)-- (5.,1.732050807568878);
\draw (5.5,0.8660254037844439)-- (5.,0.);
\draw (0.5,0.8660254037844386)-- (5.5,0.8660254037844439);
\draw (5.,1.732050807568878)-- (1.,1.7320508075688772);
\draw (2.,3.4641016151377544)-- (4.,3.4641016151377526);
\draw (4.5,2.598076211353314)-- (1.5,2.598076211353316);
\draw (2.5,4.330127018922193)-- (3.5,4.330127018922193);
\draw (3.,5.196152422706633)-- (6.,0.);
\draw (6.,0.)-- (0.,0.);
\draw (0.,0.)-- (3.,5.196152422706633);
\draw [color=aqaqaq] (0.,0.)-- (1.,0.);
\draw [color=aqaqaq] (1.,0.)-- (1.5,0.8660254037844384);
\draw [color=aqaqaq] (1.5,0.8660254037844384)-- (0.5,0.8660254037844386);
\draw [color=aqaqaq] (0.5,0.8660254037844386)-- (0.,0.);
\draw [color=aqaqaq] (6.,0.)-- (5.5,0.8660254037844439);
\draw [color=aqaqaq] (5.5,0.8660254037844439)-- (4.5,0.8660254037844394);
\draw [color=aqaqaq] (4.5,0.8660254037844394)-- (5.,0.);
\draw [color=aqaqaq] (5.,0.)-- (6.,0.);
\draw (0.5,0.8660254037844386)-- (1.,1.7320508075688772);
\draw (1.,1.7320508075688772)-- (1.5,0.8660254037844384);
\draw (1.5,0.8660254037844384)-- (0.5,0.8660254037844386);
\draw [color=aqaqaq] (3.,5.196152422706633)-- (2.5,4.330127018922193);
\draw [color=aqaqaq] (2.5,4.330127018922193)-- (3.,3.4641016151377544);
\draw [color=aqaqaq] (3.,3.4641016151377544)-- (3.5,4.330127018922193);
\draw [color=aqaqaq] (3.5,4.330127018922193)-- (3.,5.196152422706633);
\draw [color=aqaqaq] (3.5,4.330127018922193)-- (4.,3.4641016151377526);
\draw [color=aqaqaq] (4.,3.4641016151377526)-- (3.5,2.598076211353316);
\draw [color=aqaqaq] (3.5,2.598076211353316)-- (3.,3.4641016151377544);
\draw [color=aqaqaq] (3.,3.4641016151377544)-- (3.5,4.330127018922193);
\draw (4.,3.4641016151377526)-- (3.5,2.598076211353316);
\draw (3.5,2.598076211353316)-- (4.5,2.598076211353314);
\draw (4.5,2.598076211353314)-- (4.,3.4641016151377526);
\draw [color=aqaqaq] (2.,3.4641016151377544)-- (1.5,2.598076211353316);
\draw [color=aqaqaq] (1.5,2.598076211353316)-- (2.,1.732050807568877);
\draw [color=aqaqaq] (2.,1.732050807568877)-- (2.5,2.598076211353316);
\draw [color=aqaqaq] (2.5,2.598076211353316)-- (2.,3.4641016151377544);
\draw (4.5,0.8660254037844394)-- (4.,0.);
\draw (4.,0.)-- (5.,0.);
\draw (5.,0.)-- (4.5,0.8660254037844394);
\draw [color=ttffcc] (1.,1.7320508075688772)-- (3.,1.7320508075688774);
\draw [color=ttffcc] (3.,1.7320508075688774)-- (2.,0.);
\draw [color=ttffcc] (2.,0.)-- (1.,1.7320508075688772);
\draw (2.5,0.8660254037844379)-- (2.,0.);
\draw (2.,0.)-- (3.,0.);
\draw (3.,0.)-- (2.5,0.8660254037844379);
\draw (2.5,2.598076211353316)-- (3.,1.7320508075688774);
\draw (3.,1.7320508075688774)-- (2.,1.732050807568877);
\draw (2.,1.732050807568877)-- (2.5,2.598076211353316);
\draw [color=aqaqaq] (2.5,2.598076211353316)-- (3.5,2.598076211353316);
\draw [color=aqaqaq] (3.5,2.598076211353316)-- (4.,1.732050807568878);
\draw [color=aqaqaq] (4.,1.732050807568878)-- (3.,1.7320508075688774);
\draw [color=aqaqaq] (3.,1.7320508075688774)-- (2.5,2.598076211353316);
\draw (3.5,2.598076211353316)-- (4.,1.732050807568878);
\draw (4.,1.732050807568878)-- (4.5,2.598076211353314);
\draw (4.5,2.598076211353314)-- (3.5,2.598076211353316);
\draw (4.,1.732050807568878)-- (5.,1.732050807568878);
\draw (5.,1.732050807568878)-- (4.5,2.598076211353314);
\draw (4.5,2.598076211353314)-- (4.,1.732050807568878);
\draw [color=aqaqaq] (4.,1.732050807568878)-- (3.5,0.8660254037844379);
\draw [color=aqaqaq] (3.5,0.8660254037844379)-- (4.5,0.8660254037844394);
\draw [color=aqaqaq] (4.5,0.8660254037844394)-- (5.,1.732050807568878);
\draw [color=aqaqaq] (5.,1.732050807568878)-- (4.,1.732050807568878);
\draw (3.5,0.8660254037844379)-- (4.,0.);
\draw (4.,0.)-- (4.5,0.8660254037844394);
\draw (4.5,0.8660254037844394)-- (3.5,0.8660254037844379);
\draw [color=aqaqaq] (3.5,0.8660254037844379)-- (2.5,0.8660254037844379);
\draw [color=aqaqaq] (2.5,0.8660254037844379)-- (3.,0.);
\draw [color=aqaqaq] (3.,0.)-- (4.,0.);
\draw [color=aqaqaq] (4.,0.)-- (3.5,0.8660254037844379);
\draw (0.,0.)-- (3.,5.196152422706633);
\draw (3.5,4.330127018922193)-- (1.,0.);
\draw (2.,0.)-- (4.,3.4641016151377526);
\draw (4.5,2.598076211353314)-- (3.,0.);
\draw (4.,0.)-- (5.,1.732050807568878);
\draw (5.5,0.8660254037844439)-- (5.,0.);
\draw (0.5,0.8660254037844386)-- (1.,0.);
\draw (2.,0.)-- (1.,1.7320508075688772);
\draw (1.5,2.598076211353316)-- (3.,0.);
\draw (4.,0.)-- (2.,3.4641016151377544);
\draw (2.5,4.330127018922193)-- (5.,0.);
\draw (6.,0.)-- (3.,5.196152422706633);
\draw (0.,0.)-- (6.,0.);
\draw (5.5,0.8660254037844439)-- (0.5,0.8660254037844386);
\draw (1.,1.7320508075688772)-- (5.,1.732050807568878);
\draw (4.5,2.598076211353314)-- (1.5,2.598076211353316);
\draw (2.,3.4641016151377544)-- (4.,3.4641016151377526);
\draw (3.5,4.330127018922193)-- (2.5,4.330127018922193);
\begin{scriptsize}
\draw [fill=black] (0.,0.) circle (2.5pt);
\draw [fill=black] (6.,0.) circle (2.5pt);
\draw [fill=uuuuuu] (3.,5.196152422706633) circle (2.5pt);
\draw [fill=black] (1.,0.) circle (2.5pt);
\draw [fill=black] (2.,0.) circle (2.5pt);
\draw [fill=black] (3.,0.) circle (2.5pt);
\draw [fill=black] (4.,0.) circle (2.5pt);
\draw [fill=black] (5.,0.) circle (2.5pt);
\draw [fill=black] (3.,5.196152422706632) circle (2.5pt);
\draw [fill=black] (2.5,4.330127018922193) circle (2.5pt);
\draw [fill=black] (2.,3.4641016151377544) circle (2.5pt);
\draw [fill=black] (1.5,2.598076211353316) circle (2.5pt);
\draw [fill=black] (1.,1.7320508075688772) circle (2.5pt);
\draw [fill=black] (0.5,0.8660254037844386) circle (2.5pt);
\draw [fill=black] (3.5,4.330127018922193) circle (2.5pt);
\draw [fill=black] (3.,3.4641016151377544) circle (2.5pt);
\draw [fill=black] (2.5,2.598076211353316) circle (2.5pt);
\draw [fill=black] (2.,1.732050807568877) circle (2.5pt);
\draw [fill=black] (1.5,0.8660254037844384) circle (2.5pt);
\draw [fill=black] (2.5,0.8660254037844379) circle (2.5pt);
\draw [fill=black] (3.,1.7320508075688774) circle (2.5pt);
\draw [fill=black] (3.5,2.598076211353316) circle (2.5pt);
\draw [fill=black] (4.,3.4641016151377526) circle (2.5pt);
\draw [fill=black] (4.5,2.598076211353314) circle (2.5pt);
\draw [fill=black] (4.,1.732050807568878) circle (2.5pt);
\draw [fill=black] (5.,1.732050807568878) circle (2.5pt);
\draw [fill=black] (3.5,0.8660254037844379) circle (2.5pt);
\draw [fill=black] (4.5,0.8660254037844394) circle (2.5pt);
\draw [fill=black] (5.5,0.8660254037844439) circle (2.5pt);
\end{scriptsize}
\end{tikzpicture}}
\raisebox{3em}{\text{and}}
\Scale[0.5]{
\begin{tikzpicture}[line cap=round,line join=round,>=triangle 45,x=1.0cm,y=1.0cm]
\clip(-1.,-1.) rectangle (7.,6.);
\fill[color=ffffff,fill=ffffff,fill opacity=0.1] (0.,0.) -- (6.,0.) -- (3.,5.196152422706633) -- cycle;
\fill[color=aqaqaq,fill=aqaqaq,fill opacity=1.0] (0.,0.) -- (1.,0.) -- (1.5,0.8660254037844384) -- (0.5,0.8660254037844386) -- cycle;
\fill[color=aqaqaq,fill=aqaqaq,fill opacity=1.0] (6.,0.) -- (5.5,0.8660254037844439) -- (4.5,0.8660254037844394) -- (5.,0.) -- cycle;
\fill[fill=black,fill opacity=1.0] (0.5,0.8660254037844386) -- (1.,1.7320508075688772) -- (1.5,0.8660254037844384) -- cycle;
\fill[color=aqaqaq,fill=aqaqaq,fill opacity=1.0] (3.,5.196152422706633) -- (2.5,4.330127018922193) -- (3.,3.4641016151377544) -- (3.5,4.330127018922193) -- cycle;
\fill[color=aqaqaq,fill=aqaqaq,fill opacity=1.0] (3.5,4.330127018922193) -- (4.,3.4641016151377526) -- (3.5,2.598076211353316) -- (3.,3.4641016151377544) -- cycle;
\fill[fill=black,fill opacity=1.0] (4.,3.4641016151377526) -- (3.5,2.598076211353316) -- (4.5,2.598076211353314) -- cycle;
\fill[color=aqaqaq,fill=aqaqaq,fill opacity=1.0] (2.,3.4641016151377544) -- (1.5,2.598076211353316) -- (2.,1.732050807568877) -- (2.5,2.598076211353316) -- cycle;
\fill[fill=black,fill opacity=1.0] (4.5,0.8660254037844394) -- (4.,0.) -- (5.,0.) -- cycle;
\fill[color=aqaqaq,fill=aqaqaq,fill opacity=1.0] (1.,1.7320508075688772) -- (2.,1.732050807568877) -- (2.5,0.8660254037844379) -- (1.5,0.8660254037844384) -- cycle;
\fill[color=aqaqaq,fill=aqaqaq,fill opacity=1.0] (2.,0.) -- (2.5,0.8660254037844379) -- (3.5,0.8660254037844379) -- (3.,0.) -- cycle;
\fill[fill=black,fill opacity=1.0] (2.,1.732050807568877) -- (3.,1.7320508075688774) -- (2.5,0.8660254037844379) -- cycle;
\fill[fill=black,fill opacity=1.0] (2.,1.732050807568877) -- (2.5,2.598076211353316) -- (3.,1.7320508075688774) -- cycle;
\fill[fill=black,fill opacity=1.0] (3.,1.7320508075688774) -- (3.5,0.8660254037844379) -- (2.5,0.8660254037844379) -- cycle;
\fill[color=ttffcc,fill=ttffcc,fill opacity=1.0] (3.,1.7320508075688774) -- (5.,1.732050807568878) -- (4.,0.) -- cycle;
\fill[fill=black,fill opacity=1.0] (4.5,2.598076211353314) -- (4.,1.732050807568878) -- (5.,1.732050807568878) -- cycle;
\fill[color=aqaqaq,fill=aqaqaq,fill opacity=1.0] (2.5,2.598076211353316) -- (3.5,2.598076211353316) -- (4.,1.732050807568878) -- (3.,1.7320508075688774) -- cycle;
\fill[fill=black,fill opacity=1.0] (3.5,2.598076211353316) -- (4.5,2.598076211353314) -- (4.,1.732050807568878) -- cycle;
\draw (0.,0.)-- (6.,0.);
\draw [color=ffffff] (0.,0.)-- (6.,0.);
\draw [color=ffffff] (6.,0.)-- (3.,5.196152422706633);
\draw [color=ffffff] (3.,5.196152422706633)-- (0.,0.);
\draw (0.5,0.8660254037844386)-- (1.,0.);
\draw (1.,1.7320508075688772)-- (2.,0.);
\draw (1.5,2.598076211353316)-- (3.,0.);
\draw (2.,3.4641016151377544)-- (4.,0.);
\draw (5.,0.)-- (2.5,4.330127018922193);
\draw (3.5,4.330127018922193)-- (1.,0.);
\draw (2.,0.)-- (4.,3.4641016151377526);
\draw (4.5,2.598076211353314)-- (3.,0.);
\draw (4.,0.)-- (5.,1.732050807568878);
\draw (5.5,0.8660254037844439)-- (5.,0.);
\draw (0.5,0.8660254037844386)-- (5.5,0.8660254037844439);
\draw (5.,1.732050807568878)-- (1.,1.7320508075688772);
\draw (2.,3.4641016151377544)-- (4.,3.4641016151377526);
\draw (4.5,2.598076211353314)-- (1.5,2.598076211353316);
\draw (2.5,4.330127018922193)-- (3.5,4.330127018922193);
\draw (3.,5.196152422706633)-- (6.,0.);
\draw (6.,0.)-- (0.,0.);
\draw (0.,0.)-- (3.,5.196152422706633);
\draw [color=aqaqaq] (0.,0.)-- (1.,0.);
\draw [color=aqaqaq] (1.,0.)-- (1.5,0.8660254037844384);
\draw [color=aqaqaq] (1.5,0.8660254037844384)-- (0.5,0.8660254037844386);
\draw [color=aqaqaq] (0.5,0.8660254037844386)-- (0.,0.);
\draw [color=aqaqaq] (6.,0.)-- (5.5,0.8660254037844439);
\draw [color=aqaqaq] (5.5,0.8660254037844439)-- (4.5,0.8660254037844394);
\draw [color=aqaqaq] (4.5,0.8660254037844394)-- (5.,0.);
\draw [color=aqaqaq] (5.,0.)-- (6.,0.);
\draw (0.5,0.8660254037844386)-- (1.,1.7320508075688772);
\draw (1.,1.7320508075688772)-- (1.5,0.8660254037844384);
\draw (1.5,0.8660254037844384)-- (0.5,0.8660254037844386);
\draw [color=aqaqaq] (3.,5.196152422706633)-- (2.5,4.330127018922193);
\draw [color=aqaqaq] (2.5,4.330127018922193)-- (3.,3.4641016151377544);
\draw [color=aqaqaq] (3.,3.4641016151377544)-- (3.5,4.330127018922193);
\draw [color=aqaqaq] (3.5,4.330127018922193)-- (3.,5.196152422706633);
\draw [color=aqaqaq] (3.5,4.330127018922193)-- (4.,3.4641016151377526);
\draw [color=aqaqaq] (4.,3.4641016151377526)-- (3.5,2.598076211353316);
\draw [color=aqaqaq] (3.5,2.598076211353316)-- (3.,3.4641016151377544);
\draw [color=aqaqaq] (3.,3.4641016151377544)-- (3.5,4.330127018922193);
\draw (4.,3.4641016151377526)-- (3.5,2.598076211353316);
\draw (3.5,2.598076211353316)-- (4.5,2.598076211353314);
\draw (4.5,2.598076211353314)-- (4.,3.4641016151377526);
\draw [color=aqaqaq] (2.,3.4641016151377544)-- (1.5,2.598076211353316);
\draw [color=aqaqaq] (1.5,2.598076211353316)-- (2.,1.732050807568877);
\draw [color=aqaqaq] (2.,1.732050807568877)-- (2.5,2.598076211353316);
\draw [color=aqaqaq] (2.5,2.598076211353316)-- (2.,3.4641016151377544);
\draw (4.5,0.8660254037844394)-- (4.,0.);
\draw (4.,0.)-- (5.,0.);
\draw (5.,0.)-- (4.5,0.8660254037844394);
\draw [color=aqaqaq] (1.,1.7320508075688772)-- (2.,1.732050807568877);
\draw [color=aqaqaq] (2.,1.732050807568877)-- (2.5,0.8660254037844379);
\draw [color=aqaqaq] (2.5,0.8660254037844379)-- (1.5,0.8660254037844384);
\draw [color=aqaqaq] (1.5,0.8660254037844384)-- (1.,1.7320508075688772);
\draw [color=aqaqaq] (2.,0.)-- (2.5,0.8660254037844379);
\draw [color=aqaqaq] (2.5,0.8660254037844379)-- (3.5,0.8660254037844379);
\draw [color=aqaqaq] (3.5,0.8660254037844379)-- (3.,0.);
\draw [color=aqaqaq] (3.,0.)-- (2.,0.);
\draw (2.,1.732050807568877)-- (3.,1.7320508075688774);
\draw (3.,1.7320508075688774)-- (2.5,0.8660254037844379);
\draw (2.5,0.8660254037844379)-- (2.,1.732050807568877);
\draw (2.,1.732050807568877)-- (2.5,2.598076211353316);
\draw (2.5,2.598076211353316)-- (3.,1.7320508075688774);
\draw (3.,1.7320508075688774)-- (2.,1.732050807568877);
\draw (3.,1.7320508075688774)-- (3.5,0.8660254037844379);
\draw (3.5,0.8660254037844379)-- (2.5,0.8660254037844379);
\draw (2.5,0.8660254037844379)-- (3.,1.7320508075688774);
\draw [color=ttffcc] (3.,1.7320508075688774)-- (5.,1.732050807568878);
\draw [color=ttffcc] (5.,1.732050807568878)-- (4.,0.);
\draw [color=ttffcc] (4.,0.)-- (3.,1.7320508075688774);
\draw (4.5,2.598076211353314)-- (4.,1.732050807568878);
\draw (4.,1.732050807568878)-- (5.,1.732050807568878);
\draw (5.,1.732050807568878)-- (4.5,2.598076211353314);
\draw [color=aqaqaq] (2.5,2.598076211353316)-- (3.5,2.598076211353316);
\draw [color=aqaqaq] (3.5,2.598076211353316)-- (4.,1.732050807568878);
\draw [color=aqaqaq] (4.,1.732050807568878)-- (3.,1.7320508075688774);
\draw [color=aqaqaq] (3.,1.7320508075688774)-- (2.5,2.598076211353316);
\draw (3.5,2.598076211353316)-- (4.5,2.598076211353314);
\draw (4.5,2.598076211353314)-- (4.,1.732050807568878);
\draw (4.,1.732050807568878)-- (3.5,2.598076211353316);
\draw (0.,0.)-- (3.,5.196152422706633);
\draw (3.5,4.330127018922193)-- (1.,0.);
\draw (2.,0.)-- (4.,3.4641016151377526);
\draw (4.5,2.598076211353314)-- (3.,0.);
\draw (4.,0.)-- (5.,1.732050807568878);
\draw (5.5,0.8660254037844439)-- (5.,0.);
\draw (0.5,0.8660254037844386)-- (1.,0.);
\draw (2.,0.)-- (1.,1.7320508075688772);
\draw (1.5,2.598076211353316)-- (3.,0.);
\draw (4.,0.)-- (2.,3.4641016151377544);
\draw (2.5,4.330127018922193)-- (5.,0.);
\draw (6.,0.)-- (3.,5.196152422706633);
\draw (2.5,4.330127018922193)-- (3.5,4.330127018922193);
\draw (4.,3.4641016151377526)-- (2.,3.4641016151377544);
\draw (1.5,2.598076211353316)-- (4.5,2.598076211353314);
\draw (5.,1.732050807568878)-- (1.,1.7320508075688772);
\draw (0.5,0.8660254037844386)-- (5.5,0.8660254037844439);
\draw (6.,0.)-- (0.,0.);
\begin{scriptsize}
\draw [fill=black] (0.,0.) circle (2.5pt);
\draw [fill=black] (6.,0.) circle (2.5pt);
\draw [fill=uuuuuu] (3.,5.196152422706633) circle (2.5pt);
\draw [fill=black] (1.,0.) circle (2.5pt);
\draw [fill=black] (2.,0.) circle (2.5pt);
\draw [fill=black] (3.,0.) circle (2.5pt);
\draw [fill=black] (4.,0.) circle (2.5pt);
\draw [fill=black] (5.,0.) circle (2.5pt);
\draw [fill=black] (3.,5.196152422706632) circle (2.5pt);
\draw [fill=black] (2.5,4.330127018922193) circle (2.5pt);
\draw [fill=black] (2.,3.4641016151377544) circle (2.5pt);
\draw [fill=black] (1.5,2.598076211353316) circle (2.5pt);
\draw [fill=black] (1.,1.7320508075688772) circle (2.5pt);
\draw [fill=black] (0.5,0.8660254037844386) circle (2.5pt);
\draw [fill=black] (3.5,4.330127018922193) circle (2.5pt);
\draw [fill=black] (3.,3.4641016151377544) circle (2.5pt);
\draw [fill=black] (2.5,2.598076211353316) circle (2.5pt);
\draw [fill=black] (2.,1.732050807568877) circle (2.5pt);
\draw [fill=black] (1.5,0.8660254037844384) circle (2.5pt);
\draw [fill=black] (2.5,0.8660254037844379) circle (2.5pt);
\draw [fill=black] (3.,1.7320508075688774) circle (2.5pt);
\draw [fill=black] (3.5,2.598076211353316) circle (2.5pt);
\draw [fill=black] (4.,3.4641016151377526) circle (2.5pt);
\draw [fill=black] (4.5,2.598076211353314) circle (2.5pt);
\draw [fill=black] (4.,1.732050807568878) circle (2.5pt);
\draw [fill=black] (5.,1.732050807568878) circle (2.5pt);
\draw [fill=black] (3.5,0.8660254037844379) circle (2.5pt);
\draw [fill=black] (4.5,0.8660254037844394) circle (2.5pt);
\draw [fill=black] (5.5,0.8660254037844439) circle (2.5pt);
\end{scriptsize}
\end{tikzpicture}}\]
respectively, where we use the color-coding of puzzle pieces described above. It is straightforward to check that these are the only $K$-puzzles in the sense of \cite[$\mathsection 3.3$]{Vakil} for this structure constant.
\qed
\end{example}

\section{Infusion, Bender-Knuth involutions and the genomic Schur function}\label{sec:genomic_Schur}

We first define {\bf genomic infusion}. Let $T\in {\tt Gen}(\alpha)$ and
$U\in {\tt Gen}(\beta/\alpha)$ where $\alpha$ is possibly a skew shape. We think of a layered tableau $(T,U)$ that is the union of $T$ and $U$. For convenience,
the labels of $T$ will be circled. Then 
\[{\tt geninf}(T,U)=(U^\star,T^\star)\] 
is obtained by the following procedure.
Consider the largest gene $\circled \GG$ (under the $<$ order) that appears in $T$. The boxes of this gene are inner corners $I$ with respect to $U$. Now apply ${\tt jdt}_I(U)$, leaving some outer corners of $\beta$. Place into these outer corners $\circled \GG$. Now consider the second largest gene $\circled \GG'$ that appears in $T$. These will form inner corners $I'$ with respect to $U':={\tt jdt}_I(U)$. Now apply ${\tt jdt}_{I'}(U')$ again leaving some outer corners of which we will fill with $\circled \GG'$. We continue in this manner until we have exhausted all genes of $T$. The ``inner'' tableau of uncircled genes is $U^\star$ and the ``outer'' tableau of circled genes is $T^\star$.  Clearly, if $\alpha$ is a straight shape, then $U^\star$ is a genomic rectification of $U$ where the order of rectification is imposed by $T$. Furthermore: 

\begin{proposition}
\label{prop:infusioninvol}
Genomic infusion is an involution, i.e.,
\[{\tt geninf}(U^\star,T^\star)=(T,U).\] 
\end{proposition}
\begin{proof}
This follows from the fact that \emph{$K$-infusion} as defined in
\cite[$\mathsection$3.1]{Thomas.Yong:V} is an involution \cite[Theorem~3.1]{Thomas.Yong:V}, combined with 
Lemma~\ref{lem:genomic tableau_jdt}.
\end{proof}

Next we define {\bf genomic Bender-Knuth} involutions. Given a genomic tableau $V$ consider the genomic subtableau $T$ consisting of genes of family $i$ and consider the genomic subtableau $U$ consisting of genes of family $i+1$. Now define ${\tt genBK}_i(V)$ to be obtained by replacing inside $V$ the subtableau $(T,U)$ with $(U^\star,T^\star)$,
switching the labels $i$ and $i+1$, keeping all other boxes of $V$ the same (and removing any circlings). 

\begin{proposition}
\label{prop:benderknuth}
${\tt genBK}_i$ is an involution. Moreover, 
${\tt genBK}_i$ defines a bijection from the set of genomic tableaux of a shape $\nu/\lambda$ of content $\gamma=(\gamma_1,\ldots,\gamma_i,\gamma_{i+1},\ldots)$ to the set of 
genomic tableaux of shape $\nu/\lambda$ of content $\gamma=(\gamma_1,\ldots,\gamma_{i+1},\gamma_{i},\ldots)$.
\end{proposition}
\begin{proof}
The first sentence is immediate from Proposition~\ref{prop:infusioninvol}. The second sentence follows from the definition of ${\tt genBK}_i$ and the first sentence.
\end{proof}

From these genomic Bender-Knuth involutions, one can define genomic versions of M.-P.~Sch\"utzenberger's \emph{promotion} and \emph{evacuation} operators. (For the classical theory, see \cite{BPS}, specifically Theorems~2.2 and~2.9, as well as the references therein.) We do not analyze these notions further in this paper.

We explore the {\bf genomic Schur function}, which we define as
\[
U_{\nu/\lambda} := \sum_{T \in {\rm Gen}(\nu/\lambda)} {\mathbf x}^T\]
where
\[\mathbf{x}^T := \prod_i x_i^{\# \text{ genes of family $i$ in $T$}}.\]

\begin{example}\label{ex:gsf}
The polynomial $U_{31}(x_1, x_2)$ is computed by the tableaux
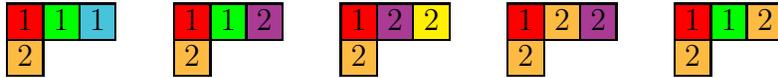
\ytableausetup{boxsize=1.1em}
\[ \ytableaushort{{*(red) 1} {*(green) 1} {*(SkyBlue) 1},{*(Dandelion) 2}} \hspace{.3in} \ytableaushort{ {*(red) 1} {*(green) 1} {*(Mulberry) 2},{*(Dandelion) 2}} \hspace{.3in} \ytableaushort{{*(red) 1} {*(Mulberry) 2} {*(Yellow) 2},{*(Dandelion) 2}}  \hspace{.3in}
\ytableaushort{{*(red) 1} {*(Dandelion) 2} {*(Mulberry) 2},{*(Dandelion) 2}} \hspace{.3in} \ytableaushort{{*(red) 1} {*(green) 1} {*(Dandelion) 2},{*(Dandelion) 2}}\]
Hence  $U_{31}(x_1, x_2) = x_1^3x_2 + x_1^2x_2^2 + x_1x_2^3 + x_1x_2^2 + x_1^2x_2 
= s_{31}(x_1,x_2) + s_{21}(x_1,x_2)$. \qed
\end{example}

\begin{theorem}\label{thm:gsf}
$U_{\nu/\lambda}\in {\sf Sym}$.
\end{theorem}
\begin{proof} The argument is an extension of the combinatorial proof of symmetry of Schur functions: It follows from Proposition~\ref{prop:benderknuth}. 
\end{proof}

Since 
\[U_{\nu/\lambda} = s_{\nu /\lambda} + \text{lower degree terms},\] 
by Theorem~\ref{thm:gsf} we have that 
$\{ U_\lambda \}$, where $\lambda$ ranges over all (straight) partitions, is a basis of ${\sf Sym}$. 

\ytableausetup{boxsize=.4em}
\begin{table}[h]
\begin{tabular}{|c|ccccccccccccccccccc|}
\hline\hline
$\lambda\backslash\mu$   & $\ydiagram{1}$ & $\ydiagram{2}$ & $\ydiagram{1,1,0}$ & $\ydiagram{3}$ & $\ydiagram{2,1}$ & $\ydiagram{1,1,1}$ & $\ydiagram{3,1}$ & $\ydiagram{2,2}$ & $\ydiagram{2,1,1}$ & $\ydiagram{3,2}$ & $\ydiagram{3,1,1}$ & $\ydiagram{2,2,1}$ & $\ydiagram{3,3}$ & $\ydiagram{3,2,1}$ & $\ydiagram{2,2,2}$ & $\ydiagram{3,3,1,0}$ & $\ydiagram{3,2,2}$ & $\ydiagram{3,3,2}$ & $\ydiagram{3,3,3}$\\ \hline
$\ydiagram{1}$ & $1$ & $\grey{0}$ & $\grey{0}$ & $\grey{0}$ & $\grey{0}$ & $\grey{0}$ & $\grey{0}$ & $\grey{0}$ & $\grey{0}$ & $\grey{0}$ & $\grey{0}$ & $\grey{0}$ & $\grey{0}$ & $\grey{0}$ & $\grey{0}$ & $\grey{0}$ & $\grey{0}$ & $\grey{0}$ & $\grey{0}$\\ \hline
$\ydiagram{2}$ & $\grey{0}$ & $1$ & $\grey{0}$ & $\grey{0}$ & $\grey{0}$ & $\grey{0}$ & $\grey{0}$ & $\grey{0}$ & $\grey{0}$ & $\grey{0}$ & $\grey{0}$ & $\grey{0}$ & $\grey{0}$ & $\grey{0}$ & $\grey{0}$ & $\grey{0}$ & $\grey{0}$ & $\grey{0}$ & $\grey{0}$\\ \hline
$\ydiagram{1,1,0}$ & $\grey{0}$ & $\grey{0}$ & $1$ & $\grey{0}$ & $\grey{0}$ & $\grey{0}$ & $\grey{0}$ & $\grey{0}$ & $\grey{0}$ & $\grey{0}$ & $\grey{0}$ & $\grey{0}$ & $\grey{0}$ & $\grey{0}$ & $\grey{0}$ & $\grey{0}$ & $\grey{0}$ & $\grey{0}$ & $\grey{0}$\\ \hline
$\ydiagram{3}$ & $\grey{0}$ & $\grey{0}$ & $\grey{0}$ & $1$ & $\grey{0}$ & $\grey{0}$ & $\grey{0}$ & $\grey{0}$ & $\grey{0}$ & $\grey{0}$ & $\grey{0}$ & $\grey{0}$ & $\grey{0}$ & $\grey{0}$ & $\grey{0}$ & $\grey{0}$ & $\grey{0}$ & $\grey{0}$ & $\grey{0}$\\ \hline
$\ydiagram{2,1,0}$ & $\grey{0}$ & $\grey{0}$ & $1$ & $\grey{0}$ & $1$ & $\grey{0}$ & $\grey{0}$ & $\grey{0}$ & $\grey{0}$ & $\grey{0}$ & $\grey{0}$ & $\grey{0}$ & $\grey{0}$ & $\grey{0}$ & $\grey{0}$ & $\grey{0}$ & $\grey{0}$ & $\grey{0}$ & $\grey{0}$\\ \hline
$\ydiagram{1,1,1,0}$ & $\grey{0}$ & $\grey{0}$ & $\grey{0}$ & $\grey{0}$ & $\grey{0}$ & $1$ & $\grey{0}$ & $\grey{0}$ & $\grey{0}$ & $\grey{0}$ & $\grey{0}$ & $\grey{0}$ & $\grey{0}$ & $\grey{0}$ & $\grey{0}$ & $\grey{0}$ & $\grey{0}$ & $\grey{0}$ & $\grey{0}$\\ \hline
$\ydiagram{3,1,0}$ & $\grey{0}$ & $\grey{0}$ & $\grey{0}$ & $\grey{0}$ & $1$ & $\grey{0}$ & $1$ & $\grey{0}$ & $\grey{0}$ & $\grey{0}$ & $\grey{0}$ & $\grey{0}$ & $\grey{0}$ & $\grey{0}$ & $\grey{0}$ & $\grey{0}$ & $\grey{0}$ & $\grey{0}$ & $\grey{0}$\\ \hline
$\ydiagram{2,2,0}$ & $\grey{0}$ & $\grey{0}$ & $\grey{0}$ & $\grey{0}$ & $\grey{0}$ & $1$ & $\grey{0}$ & $1$ & $\grey{0}$ & $\grey{0}$ & $\grey{0}$ & $\grey{0}$ & $\grey{0}$ & $\grey{0}$ & $\grey{0}$ & $\grey{0}$ & $\grey{0}$ & $\grey{0}$ & $\grey{0}$\\ \hline
$\ydiagram{2,1,1,0}$ & $\grey{0}$ & $\grey{0}$ & $\grey{0}$ & $\grey{0}$ & $\grey{0}$ & $2$ & $\grey{0}$ & $\grey{0}$ & $1$ & $\grey{0}$ & $\grey{0}$ & $\grey{0}$ & $\grey{0}$ & $\grey{0}$ & $\grey{0}$ & $\grey{0}$ & $\grey{0}$ & $\grey{0}$ & $\grey{0}$\\ \hline
$\ydiagram{3,2,0}$ & $\grey{0}$ & $\grey{0}$ & $\grey{0}$ & $\grey{0}$ & $\grey{0}$ & $1$ & $\grey{0}$ & $1$ & $1$ & $1$ & $\grey{0}$ & $\grey{0}$ & $\grey{0}$ & $\grey{0}$ & $\grey{0}$ & $\grey{0}$ & $\grey{0}$ & $\grey{0}$ & $\grey{0}$\\ \hline
$\ydiagram{3,1,1,0}$ & $\grey{0}$ & $\grey{0}$ & $\grey{0}$ & $\grey{0}$ & $\grey{0}$ & $1$ & $\grey{0}$ & $\grey{0}$ & $2$ & $\grey{0}$ & $1$ & $\grey{0}$ & $\grey{0}$ & $\grey{0}$ & $\grey{0}$ & $\grey{0}$ & $\grey{0}$ & $\grey{0}$ & $\grey{0}$\\ \hline
$\ydiagram{2,2,1,0}$ & $\grey{0}$ & $\grey{0}$ & $\grey{0}$ & $\grey{0}$ & $\grey{0}$ & $1$ & $\grey{0}$ & $\grey{0}$ & $1$ & $\grey{0}$ & $\grey{0}$ & $1$ & $\grey{0}$ & $\grey{0}$ & $\grey{0}$ & $\grey{0}$ & $\grey{0}$ & $\grey{0}$ & $\grey{0}$\\ \hline
$\ydiagram{3,3,0}$ & $\grey{0}$ & $\grey{0}$ & $\grey{0}$ & $\grey{0}$ & $\grey{0}$ & $\grey{0}$ & $\grey{0}$ & $\grey{0}$ & $\grey{0}$ & $\grey{0}$ & $\grey{0}$ & $1$ & $1$ & $\grey{0}$ & $\grey{0}$ & $\grey{0}$ & $\grey{0}$ & $\grey{0}$ & $\grey{0}$\\ \hline
$\ydiagram{3,2,1,0}$ & $\grey{0}$ & $\grey{0}$ & $\grey{0}$ & $\grey{0}$ & $\grey{0}$ & $1$ & $\grey{0}$ & $\grey{0}$ & $2$ & $\grey{0}$ & $1$ & $2$ & $\grey{0}$ & $1$ & $\grey{0}$ & $\grey{0}$ & $\grey{0}$ & $\grey{0}$ & $\grey{0}$\\ \hline
$\ydiagram{2,2,2,0}$ & $\grey{0}$ & $\grey{0}$ & $\grey{0}$ & $\grey{0}$ & $\grey{0}$ & $\grey{0}$ & $\grey{0}$ & $\grey{0}$ & $\grey{0}$ & $\grey{0}$ & $\grey{0}$ & $\grey{0}$ & $\grey{0}$ & $\grey{0}$ & $1$ & $\grey{0}$ & $\grey{0}$ & $\grey{0}$ & $\grey{0}$\\ \hline
$\ydiagram{3,3,1,0}$ & $\grey{0}$ & $\grey{0}$ & $\grey{0}$ & $\grey{0}$ & $\grey{0}$ & $\grey{0}$ & $\grey{0}$ & $\grey{0}$ & $\grey{0}$ & $\grey{0}$ & $\grey{0}$ & $1$ & $\grey{0}$ & $1$ & $\grey{0}$ & $1$ & $\grey{0}$ & $\grey{0}$ & $\grey{0}$\\ \hline
$\ydiagram{3,2,2,0}$ & $\grey{0}$ & $\grey{0}$ & $\grey{0}$ & $\grey{0}$ & $\grey{0}$ & $\grey{0}$ & $\grey{0}$ & $\grey{0}$ & $\grey{0}$ & $\grey{0}$ & $\grey{0}$ & $\grey{0}$ & $\grey{0}$ & $\grey{0}$ & $2$ & $\grey{0}$ & $1$ & $\grey{0}$ & $\grey{0}$\\ \hline
$\ydiagram{3,3,2,0}$ & $\grey{0}$ & $\grey{0}$ & $\grey{0}$ & $\grey{0}$ & $\grey{0}$ & $\grey{0}$ & $\grey{0}$ & $\grey{0}$ & $\grey{0}$ & $\grey{0}$ & $\grey{0}$ & $\grey{0}$ & $\grey{0}$ & $\grey{0}$ & $1$ & $\grey{0}$ & $1$ & $1$ & $\grey{0}$\\ \hline
$\ydiagram{3,3,3,0}$ & $\grey{0}$ & $\grey{0}$ & $\grey{0}$ & $\grey{0}$ & $\grey{0}$ & $\grey{0}$ & $\grey{0}$ & $\grey{0}$ & $\grey{0}$ & $\grey{0}$ & $\grey{0}$ & $\grey{0}$ & $\grey{0}$ & $\grey{0}$ & $\grey{0}$ & $\grey{0}$ & $\grey{0}$ & $\grey{0}$ & $1$\\ \hline
\end{tabular}
\caption{Transition matrix between the $\{U_{\lambda}(x_1,x_2,x_3)\}$ to $\{s_{\mu}(x_1,x_2,x_3)\}$ bases.
} 
\label{tab:transition}
\end{table}

While in small examples $U_{\nu/\lambda}$ is Schur-positive (cf.\ Table~\ref{tab:transition}), this is not true in general:

\begin{example}\label{ex:gsf_not_s-positive}
One may check that $38$ tableaux contribute to $U_{333}(x_1, x_2, x_3, x_4)$. Expanding this polynomial in the Schur basis yields 
\begin{align*}
U_{333}(x_1, x_2, x_3, x_4) = s_{333}(x_1, x_2, x_3, x_4) &+ s_{3221}(x_1, x_2, x_3, x_4) \\ &+ s_{2221}(x_1, x_2, x_3, x_4) - s_{2222}(x_1, x_2, x_3, x_4). \quad \quad \ \ \  \qed
\end{align*}
\end{example}
\ytableausetup{boxsize=1.1em}

Also, the structure coefficients for the $U$-basis do not possess any positivity or alternating positivity properties:
\begin{example}
Using Table~\ref{tab:transition}, one can check that $U_{22}\cdot U_1=U_{32}+U_{221}-U_{22}-U_{111}$.\qed
\end{example}
At present, we are unaware of any geometric significance of these polynomials.

\section{Shifted genomic tableaux}\label{sec:shifted_genomic_tableaux}
\ytableausetup{aligntableaux=top}

Recall, the {\bf shifted diagram} of a strictly decreasing partition is given by taking the ordinary Young diagram and indenting row $i$ (from the top) $i-1$ positions to the right.
Let 
\[\mathcal{D}:=\{1' < 1 < 2' < 2 < \cdots\}.\] 
A {\bf $P$-tableau} is a filling of shifted shape $\nu / \lambda$ 
with entries from $\mathcal{D}$ such that:
\begin{itemize}
\item[(P.1)] rows and columns weakly increase (left to right, top to bottom);
\item[(P.2)] each unprimed letter appears at most once in any column;
\item[(P.3)]  each primed letter appears at most once in any row; and
\item[(P.4)] every primed letter $k'$ has an unprimed $k$ southwest of it. 
\end{itemize}
The \emph{Schur $P$-function} $P_\lambda$ is a generating function over these tableaux (for more history and development of these functions, see e.g., \cite{hoffman.humphreys} or\cite{Stembridge}).

\begin{example}
$\begin{ytableau}
1 & 2' & 3\\
\none & 2\\
\end{ytableau}$ is a $P$-tableau of shape $\lambda=(3,1)$. The tableau $\ytableaushort{2 {3'} 4 4, \none {3'} 6, \none \none 7}$ is not a $P$-tableau because it violates both (P.3) and (P.4). 
However, if the lower $3'$ changes to $3$, the result is a $P$-tableau. \qed
\end{example}

For $\alpha \in \mathcal{D}$, write $|\alpha| = k$ if $\alpha \in \{ k' , k \}$. We use initial letters of the Greek alphabet ($\alpha, \beta, \gamma, \ldots$) for elements of $\mathcal{D}$, reserving Roman letters for elements of $\Z$.

For fixed $k \in \Z_{\geq 0}$, place a total order $\prec$ on those boxes with entry $k'$ in top to bottom order and on those boxes with entry $k$ using left to right order; declare the boxes containing $k'$ to precede those containing $k$.
A {\bf gene} (of family $k$) in a
$P$-tableau $T$ is a set $\GG$ of boxes of $T$ such that:
\begin{itemize}
\item each entry in $\GG$ is $k'$ or $k$;
\item the boxes of $\GG$ are consecutive in the $\prec$-order; and
\item no two boxes of $\GG$ appear in the same row or the same column.
\end{itemize}
We write ${\tt family}(\GG) = k$.

\begin{example}
Consider the following three colorings of the same $P$-tableau: 
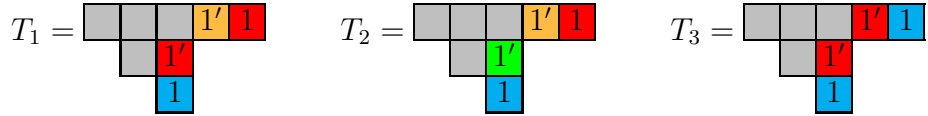
\[ T_1 = \ytableaushort{{*(lightgray)\blank} {*(lightgray)\blank} {*(lightgray)\blank}{*(Dandelion) 1'} {*(red) 1}, \none {*(lightgray)\blank} {*(red) 1'}, \none \none {*(cyan) 1} } \hspace{1cm} T_2=\ytableaushort{{*(lightgray)\blank} {*(lightgray)\blank} {*(lightgray)\blank}{*(Dandelion) 1'} {*(red) 1}, \none {*(lightgray)\blank} {*(green) 1'}, \none \none {*(cyan) 1} } \hspace{1cm} T_3=\ytableaushort{{*(lightgray)\blank} {*(lightgray)\blank} {*(lightgray)\blank}{*(red) 1'} {*(cyan) 1}, \none {*(lightgray)\blank} {*(red) 1'}, \none \none {*(cyan) 1} }\]
The {\textcolor{red}{red boxes}} in $T_1$ do not form a gene, since they are not consecutive in $\prec$-order (in view of the {\textcolor{cyan}{blue $1$}}).
In $T_2$ and $T_3$, the boxes of each color form valid genes.
\qed
\end{example}

A {\bf genomic $P$-tableau} is a $P$-tableau $T$ together with a partition of its boxes into genes such that for every primed box ${\sf b}$, there is an box ${\sf c}$ that is weakly southwest of ${\sf b}$ from a different gene than ${\sf b}$ but of the same family.
The {\bf content} of $T$ is the number of genes of each family. A {\bf genotype} $G$ of $T$ is a choice of a single box from each gene. Depict $G$ by erasing the entries in boxes that are not chosen. A $P$-tableau may be identified with the genomic $P$-tableau where each box is its own gene.

\begin{example}\label{ex:B.genomic.tableau}
Let $\nu=(6,4,1)$ and $\lambda=(4,2)$. Then a genomic $P$-tableau $T$ of shape
$\nu/\lambda$ and its two genotypes $G_1, G_2$ are
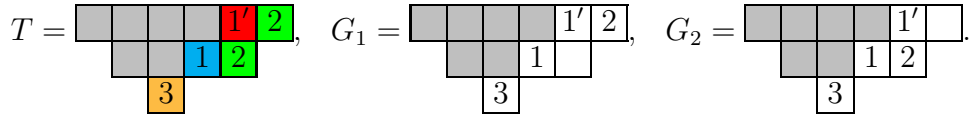
\[T=\begin{ytableau}
*(lightgray)\blank & *(lightgray)\blank & *(lightgray)\blank & *(lightgray)\blank & *(red) 1' & *(green) 2\\
\none & *(lightgray)\blank & *(lightgray)\blank & *(cyan) 1 & *(green) 2\\
\none & \none & *(Dandelion) 3\\
\end{ytableau}, \ \ \
G_1=\begin{ytableau}
*(lightgray)\blank & *(lightgray)\blank & *(lightgray)\blank & *(lightgray)\blank & 1' & 2\\
\none & *(lightgray)\blank & *(lightgray)\blank & 1 & \blank\\
\none & \none & 3\\
\end{ytableau}, \ \ \
G_2=\begin{ytableau}
*(lightgray)\blank & *(lightgray)\blank & *(lightgray)\blank & *(lightgray)\blank & 1' & \blank\\
\none & *(lightgray)\blank & *(lightgray)\blank & 1 & 2\\
\none & \none & 3\\
\end{ytableau}.
\]
The content of $T$ is $\mu =(2,1,1)$. \qed
\end{example}

Given a word $w$ using the alphabet $\mathcal{D}$, $\hat{w}$ is the
word obtained by writing $w$ backwards, and
replacing each $k'$ with $k$ while simultaneously replacing
each $k$ with $(k + 1)'$. Let 
\[{\tt doubleseq}(G) :={\tt seq}(G){\widehat{{\tt seq}(G)}}.\]
Say ${\tt doubleseq}(G)$ is {\bf locally ballot} at the letter $\alpha \in \mathcal{D}$, 
if $|\alpha| = 1$ or if in ${\tt doubleseq}(G)$ the number of $|\alpha|$'s appearing strictly before that $\alpha$ 
is strictly less than the number of $(|\alpha| - 1)$'s appearing strictly before that $\alpha$. 
Declare ${\tt doubleseq}(G)$ to be {\bf ballot} if 
it is locally ballot at each letter. Finally, $G$ is {\bf ballot} if ${\tt doubleseq}(G)$ is ballot, and the genomic $P$-tableau $T$ is {\bf ballot} if every genotype of $T$ is ballot. 

\begin{example}
Let $G_1$ and $G_2$ be as in Example~\ref{ex:B.genomic.tableau}. Then
\[{\tt doubleseq}(G_1)=2 1' 1 3 4' 2' 1 3' \text{\ and ${\tt doubleseq}(G_2)=1' 2 1 3 4'2' 3' 1$}.\] 
The former is not ballot, as it starts with $2$. Hence the genomic $P$-tableau $T$ of Example~\ref{ex:B.genomic.tableau} 
is not ballot. $G_2$ is also not ballot: ${\tt doubleseq}(G_2)$ 
is locally ballot at every position except the $2$ in second position; although there is a $1'$ before this $2$, there is no $1$. To emphasize the differences between ballotness in this section versus ballotness in Section~\ref{sec:sequences}, note that deleting the primes gives $12134231$, which \emph{is} ballot in the earlier sense.
\qed
\end{example}

A {\bf $Q$-tableau} is a filling of $\nu / \lambda$ with entries from $\mathcal{D}$ satisfying (P.1)--(P.3) and
\begin{itemize}
\item[(Q.4)] no primed letters appear on the main diagonal.
\end{itemize}
(Observe that (Q.4) is a weakening of (P.4), so a $P$-tableau is a $Q$-tableau.)

A {\bf gene} (of family $k$) in a $Q$-tableau $T$ is a set $\GG$ of boxes such that:
\begin{itemize}
\item each entry of $\GG$ is $k'$ or $k$,
\item the boxes of $\GG$ are consecutive in the $\prec$-order, and 
\item no two boxes of $\GG$ \emph{with the same label} appear in the same row or the same column. 
\end{itemize}

We write ${\tt family}(\GG) = k$ as before.

A {\bf genomic $Q$-tableau} is a  $Q$-tableau $T$ together with a partition of its boxes into genes. The definition of {\bf ballotness} for genomic $Q$-tableaux is the same as for genomic $P$-tableaux. Let $\Pgenomictableau{\nu/\lambda}{\mu}$ and $\Qgenomictableau{\nu/\lambda}{\mu}$ respectively denote the sets of genomic $P$- and $Q$-tableaux of shape $\nu/\lambda$ and content $\mu$.

\begin{lemma}
$\Pgenomictableau{\nu/\lambda}{\mu} \subseteq \Qgenomictableau{\nu/\lambda}{\mu}$.
\end{lemma}
\begin{proof}
Let $T \in \Pgenomictableau{\nu/\lambda}{\mu}$.
The definition of a gene in a $Q$-tableau differs from that for $P$-tableaux only in that it allows $k'$ and $k$ in the same row or column to be in the same gene. Hence each gene of $T$ is a gene in the $Q$-tableau sense. Thus $T \in \Qgenomictableau{\nu/\lambda}{\mu}$.
\end{proof}

In the announcement version of this paper (\cite[$\mathsection$3.2]{PY:FPSAC}), we used tableaux of  
self-conjugate shape rather than shifted shape. This is equivalent
to our present discussion, as those tableaux may be recovered by reflecting the shifted tableaux.

\section{Maximal orthogonal and Lagrangian Grassmannians}
\label{sec:result_B}

Let ${\sf G}/{\sf P}$ be a {\bf generalized flag variety}, where
${\sf G}$ is a complex, connected, reductive Lie group and ${\sf P}$ is a parabolic subgroup 
containing a Borel subgroup 
${\sf B}$. Let ${\sf B}_{-}$ be the opposite Borel to ${\sf B}$ with respect to  a choice of maximal torus ${\sf T}\subseteq {\sf B}$. The {\bf Schubert cells} of ${\sf G}/{\sf P}$ are the ${\sf B}_{-}$-orbits, and the 
{\bf Schubert varieties} $V_{\lambda}$ are their closures. Here $\lambda\in
W/W_{\sf P}$ where $W$ is the Weyl group of ${\sf G}$
and $W_{\sf P}$ is the parabolic subgroup of $W$ corresponding to ${\sf P}$. The classes of Schubert
structure sheaves $\{[{\mathcal O}_{V_{\lambda}}]\}$ 
form a ${\mathbb Z}$-linear basis of the {\bf Grothendieck ring} $K^{0}({\sf G}/{\sf P})$. Let $t_{\lambda,\mu}^{\nu}$ be the structure constants with respect 
to this basis. A.~Buch \cite[Conjecture~9.2]{Buch:KLR} conjectured the sign-alternation:
\[(-1)^{{\rm codim}_{{\sf G}/{\sf P}}(V_\nu)-
{\rm codim}_{{\sf G}/{\sf P}}(V_\lambda)-
{\rm codim}_{{\sf G}/{\sf P}}(V_\mu)}t_{\lambda,\mu}^{\nu}\geq 0.\] 
This was subsequently
proved by M.~Brion \cite{Brion:positivity}. While the Grassmannian $X$ is the most well-studied case of ${\sf G}/{\sf P}$, we now turn to an investigation of the next two most important cases when ${\sf P}$ is maximal parabolic.

Fix a non-degenerate, symmetric bilinear form $\beta(\cdot, \cdot)$ on ${\mathbb C}^{2n+1}$.
A subspace $V\subseteq {\mathbb C}^{2n+1}$ is {\bf isotropic} with respect to  $\beta$ 
if $\beta(\vec v, \vec w) = 0$ for all $\vec v, \vec w\in V$.
Let 
\[Y={\rm OG}(n,2n+1)\] 
be the {\bf maximal orthogonal Grassmannian}, i.e., the parameter space
 of all such isotropic $n$-dimensional subspaces in ${\mathbb C}^{2n+1}$.
Define the {\bf shifted staircase} $\delta_n$ to be the shifted shape
whose $i$th row is of length $i$ for $1\leq i\leq n$.
The Schubert varieties $Y_{\lambda}$ of $Y$
are indexed by shifted Young diagrams 
\[\lambda=(\lambda_1>\lambda_2>\dots>\lambda_n)\] 
contained
in $\delta_n$,  i.e., 
\[\lambda_k\leq n-k+1 \text{\ for $1\leq k\leq n$.}\]
We have 
\[{\rm codim}_{Y}(Y_{\lambda})=|\lambda|.\]
Let $b_{\lambda,\mu}^{\nu}$ be $t_{\lambda,\mu}^{\nu}$ in this case.
The first combinatorial rule for $b_{\lambda,\mu}^{\nu}$ was conjectured
in \cite{Thomas.Yong:V} and proved in \cite{Clifford.Thomas.Yong}, using \cite{Buch.Ravikumar}.

The following is a new rule for these structure coefficients. This rule directly extends the rule of J.~Stembridge \cite[Theorem~8.3]{Stembridge} for the ordinary cohomological structure constants of $Y$. (J.~Stembridge's rule is stated in terms of projective representation theory of symmetric groups; the application to $H^\star(Y)$ is due to P.~Pragacz \cite{pragacz}.) 

\begin{theorem}[OG Genomic Littlewood-Richardson rule]\label{thm:B_lr_rule}
\[b_{\lambda, \mu}^\nu = 
(-1)^{|\nu| - |\lambda| - |\mu|} \text{times the number of  ballot genomic $P$-tableaux of shape $\nu / \lambda$ with content $\mu$}.\]
\end{theorem}

\begin{example}\label{ex:CTY}
(cf. \cite[Example 1.3]{Clifford.Thomas.Yong})
That \[b_{(3,1), (3,1)}^{(5,3,1)}({\rm OG}(n,2n+1)) = -6\] is witnessed by:

{\small \[\ytableaushort{{*(lightgray) \blank} {*(lightgray) \blank} {*(lightgray) \blank} {*(red) 1'} {*(SkyBlue) 1}, \none {*(lightgray) \blank} {*(red) 1'} {*(Dandelion) 2}, \none \none {*(green) 1}} 
\hspace{.3in} 
\ytableaushort{{*(lightgray) \blank} {*(lightgray) \blank} {*(lightgray) \blank} {*(green) 1} {*(SkyBlue) 1}, \none {*(lightgray) \blank} {*(red) 1'} {*(Dandelion) 2}, \none \none {*(green) 1}} 
\hspace{.3in}
\ytableaushort{{*(lightgray) \blank} {*(lightgray) \blank} {*(lightgray) \blank} {*(red) 1'} {*(SkyBlue) 1}, \none {*(lightgray) \blank} {*(green) 1} {*(Dandelion) 2}, \none \none {*(Dandelion) 2}}\] 
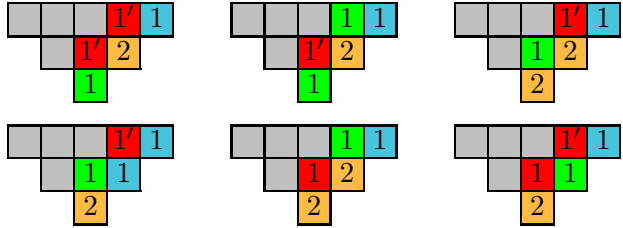
\[\ytableaushort{{*(lightgray) \blank} {*(lightgray) \blank} {*(lightgray) \blank} {*(red) 1'} {*(SkyBlue) 1}, \none {*(lightgray) \blank} {*(green) 1} {*(SkyBlue) 1}, \none \none {*(Dandelion) 2}} 
\hspace{.3in} 
\ytableaushort{{*(lightgray) \blank} {*(lightgray) \blank} {*(lightgray) \blank} {*(green) 1} {*(SkyBlue) 1}, \none {*(lightgray) \blank} {*(red) 1} {*(Dandelion) 2}, \none \none {*(Dandelion) 2}} 
\hspace{.3in} 
\ytableaushort{{*(lightgray) \blank} {*(lightgray) \blank} {*(lightgray) \blank} {*(red) 1'} {*(SkyBlue) 1}, \none {*(lightgray) \blank} {*(red) 1} {*(green) 1}, \none \none {*(Dandelion) 2}}\] }
\qed
\end{example}

Fix a symplectic bilinear form $\omega(\cdot, \cdot)$ on $\mathbb{C}^{2n}$. The {\bf Lagrangian Grassmannian} 
\[Z = {\rm LG}(n,2n)\] 
is the parameter space of $n$-dimensional linear subspaces of $\mathbb{C}^{2n}$ that are isotropic with respect to $\omega$.
The Schubert varieties $\{Z_\lambda\}$ of $Z$ are indexed by the same shifted Young diagrams $\lambda$ as above; also, ${\rm codim}_Z(Z_\lambda) = |\lambda|$. 
Let $c_{\lambda,\mu}^{\nu}$ be $t_{\lambda, \mu}^\nu$ in this case. 

There is a well-known relationship
in the ``cohomological case'', i.e., when $|\lambda|+|\mu|=|\nu|$,
between the structure constants for $Y$ and $Z$:
\begin{equation}
\label{eqn:cohomolog}
c_{\lambda,\mu}^{\nu} = 2^{\ell(\lambda) + \ell(\mu) - \ell(\nu)} b_{\lambda,\mu}^{\nu},
\end{equation} 
where $\ell(\pi)$ denotes the number of nonzero parts of $\pi$. We are not
aware of any generalization of (\ref{eqn:cohomolog}); cf.~\cite[Examples~4.9 and 5.8]{Buch.Ravikumar}. On the other hand, we propose the following extension of this
relationship:

\begin{conjecture}
\label{conj:ogvskl}
For any strict partitions $\lambda, \mu, \nu$, we have
$|b_{\lambda, \mu}^\nu| \leq |c_{\lambda, \mu}^\nu|$. \qed
\end{conjecture}

This conjecture is true in the cohomological case since it is known that
$\ell(\lambda)+\ell(\mu)\geq \ell(\nu)$ whenever $b_{\lambda,\mu}^{\nu}>0$. Moreover, we have verified this conjecture by computer for $n\leq 6$. 
In addition,
by \cite{Buch.Ravikumar}, this conjecture holds whenever $\mu$ has a single part. 

Let $\Qballot{\nu/\lambda}{\mu} := \{ \text{ballot genomic $Q$-tableaux of shape } \nu / \lambda \text{ with content } \mu \}$.

\begin{conjecture}\label{conj:upper_bound}
$|c^\nu_{\lambda,\mu}| \leq \#\Qballot{\nu/\lambda}{\mu}$.
\qed
\end{conjecture}

\begin{example}
Let $\lambda = (3,1), \mu = (2,1)$ and $\nu = (4,3,1)$. Then $\#\Qballot{\nu/\lambda}{\mu} = 6$:
{\small \[
\ytableaushort{\blank \blank \blank {*(Cyan) 1} , \none \blank {*(Red) 1'} {*(Yellow) 2} , \none \none {*(Red) 1} }
\, \,
\ytableaushort{\blank \blank \blank {*(Cyan) 1} , \none \blank {*(Red) 1'} {*(Yellow) 2'} , \none \none {*(Red) 1} }
\, \,
\ytableaushort{\blank \blank \blank {*(Cyan) 1} , \none \blank {*(Red) 1} {*(Yellow) 2} , \none \none {*(Yellow) 2} }
\, \,
\ytableaushort{\blank \blank \blank {*(Cyan) 1} , \none \blank {*(Red) 1} {*(Yellow) 2'} , \none \none {*(Yellow) 2} }
\, \,
\ytableaushort{\blank \blank \blank {*(Cyan) 1} , \none \blank {*(Red) 1'} {*(Yellow) 2} , \none \none {*(Yellow) 2} }
\, \,
\ytableaushort{\blank \blank \blank {*(Cyan) 1} , \none \blank {*(Red) 1'} {*(Yellow) 2'} , \none \none {*(Yellow) 2} }
\, \,
\] }
The third tableau above is the only one that is a genomic $P$-tableau; hence $b^\nu_{\lambda,\mu} = -1$. 
Therefore Conjectures~\ref{conj:upper_bound} and~\ref{conj:ogvskl} predict $1  \leq |c^\nu_{\lambda,\mu}| \leq 6$. Indeed, $c^\nu_{\lambda,\mu} = -5$.
\qed

\end{example}

We have computer verified Conjecture~\ref{conj:upper_bound} for $n \leq 6$. Moreover, the bound is sharp, as indicated in the
two propositions below.

\begin{proposition}
For $\mu = (p)$, $|c^\nu_{\lambda,\mu}| = \#\Qballot{\nu/\lambda}{\mu}$.
\end{proposition}
\begin{proof}
By applying $\Gamma$ (defined in Section~\ref{sec:shifted_standardization}) to the tableaux in $\Qballot{\nu/\lambda}{{(p)}}$ and retaining the primes, one obtains precisely the \emph{KLG-tableaux} of A.~Buch--V.~Ravikumar \cite[$\mathsection$5]{Buch.Ravikumar}. By \cite[Corollary~5.6]{Buch.Ravikumar}, the number of the latter is $(-1)^{|\nu|-|\lambda|-p} c^\nu_{\lambda,(p)}$.
\end{proof}

\begin{proposition}
For $|\nu| \leq |\lambda| + |\mu|$, $|c^\nu_{\lambda,\mu}| = \#\Qballot{\nu/\lambda}{\mu}$.
\end{proposition}
\begin{proof}
When $|\nu| < |\lambda| + |\mu|$, $c^\nu_{\lambda,\mu} = 0$ for geometric reasons. Clearly in this case also $\Qballot{\nu/\lambda}{\mu} = \emptyset$.

Suppose $|\nu| = |\lambda| + |\mu|$ and $T \in \Qballot{\nu / \lambda}{\mu}$. 
The number of boxes of $\nu / \lambda$ on the main diagonal is $\ell(\nu) - \ell(\lambda)$.
By pigeonhole, each gene of $T$ is a single box. Hence these tableaux are exactly the tableaux of \cite[Theorem~8.3]{Stembridge} with condition (2) removed. Therefore by the discussion of \cite[p.~126]{Stembridge}, $\#\Qballot{\nu/\lambda}{\mu}$ is the coefficient of the Schur $Q$-function $Q_\mu$ in the expansion of the skew Schur $Q$-function $Q_{\nu / \lambda}$. It is well known that these coefficients agree with the structure constants for $Z$ in this case.
\end{proof}

That is, we conjecturally have combinatorially-related upper and lower bounds for $|c_{\lambda, \mu}^\nu|$
in terms of genomic tableaux. Let \[\Pballot{\nu / \lambda}{\mu}:=\{\text{ballot genomic $P$-tableaux of} \text{ shape $\nu / \lambda$} \text{ with content }\mu\}.\] Naturally, one seeks a set 
${\tt QBallot}^\star_\mu (\nu/\lambda)$ satisfying
\[{\tt PBallot}_\mu (\nu/\lambda) \subseteq {\tt QBallot}^\star_\mu (\nu/\lambda) \subseteq {\tt QBallot}_\mu (\nu/\lambda),\] such that $\#{\tt QBallot}^\star_\mu (\nu/\lambda) = |c^\nu_{\lambda,\mu}|$.  Let 
\begin{multline}
{\tt QBallot}^\dagger_\mu (\nu/\lambda) := \\ \nonumber
\{ T \in {\tt QBallot}_\mu (\nu/\lambda) : \text{no gene contains both primed and unprimed labels} \}.
\end{multline}

\begin{conjecture}
$\#{\tt QBallot}^\dagger_\mu (\nu/\lambda) \leq |c^\nu_{\lambda,\mu}|$.
\qed
\end{conjecture}

This has also been computer-checked for $n \leq 6$. 
It suggests that one should look to define ${\tt QBallot}^\star_\mu(\nu / \lambda)$ from ${\tt QBallot}_\mu(\nu / \lambda)$ by imposing a condition on genes with both primed and unprimed labels.

\section{Proof of OG Genomic Littlewood-Richardson rule (Theorem~\ref{thm:B_lr_rule})}

Our proof of Theorem~\ref{thm:B_lr_rule} proceeds parallel to the first proof of
Theorem~\ref{thm:A_lr_rule}. (We are not aware of any set-valued tableau or puzzle formulation of Theorem~\ref{thm:B_lr_rule}.)

\subsection{Shifted $K$-(semi)standardization maps}\label{sec:shifted_standardization}
Let $T$ be a genomic $P$-tableau. Impose a total order on genes of $T$ by $\GG_1 < \GG_2$ if ${\sf b}_1 \prec {\sf b}_2$, for ${\sf b}_i$ a box of $\GG_i$. (Note that since the boxes of a gene form a $\prec$-interval, this order is well-defined.)
 
A {\bf shifted increasing tableau} is a filling of a shifted shape that strictly
increases along rows and down columns (see \cite[$\mathsection$7]{Thomas.Yong:V} and \cite{Clifford.Thomas.Yong}). Define the
{\bf shifted $K$-standardization map}
\[\Gamma:\Pgenomictableau{\nu/\lambda}{\mu}\to {\tt Inc}(\nu/\lambda)\]
by filling the $i$th gene in $<$-order with the entry $i$.

\begin{example}\label{ex:phi_type_B}
If $T$ is the genomic $P$-tableau 
$\begin{ytableau}
*(lightgray)\blank & *(lightgray)\blank & *(lightgray)\blank & *(lightgray)\blank & *(red) 1' & *(green) 2\\
\none & *(lightgray)\blank & *(lightgray)\blank & *(cyan) 1 & *(green) 2\\
\none & \none & *(Dandelion) 3\\
\end{ytableau}$
in Example~\ref{ex:B.genomic.tableau}, then 
\[\Gamma(T) = \ytableaushort{{*(lightgray)\blank} {*(lightgray)\blank} {*(lightgray)\blank} {*(lightgray)\blank} 1 3, \none {*(lightgray)\blank} {*(lightgray)\blank} 2 3, \none \none 4}.\]
\qed
\end{example}

Recall
\[
{\mathcal P}_k(\mu):=\left\{ 1 + \sum_{i<k} \mu_i, 2 + \sum_{i<k} \mu_i, \dots, \sum_{j\leq k} \mu_j\right\}.
\]
and let $S\in {\tt Inc}(\nu/\lambda)$ have largest entry $n$. Let 
\[\mu = (\mu_1, \mu_2, \dots, \mu_h)\] 
be a composition of $n$.
The {\bf shifted $K$-semistandardization} $\Delta_\mu(S)$ with respect to $\mu$ is defined as follows.
Replace each entry $i$ in $S$ with $k_i$ for the unique $k$ such that $i \in \mathcal{P}_k(\mu)$.
For each $k_h$, replace it with $k'$ if there is a $k_j$ southwest of it with $h < j$; otherwise replace it with $k$. If the result is a $P$-tableau, define a 
(putative) genomic $P$-tableau structure by putting all boxes that have the same entry in $S$ into the same gene. If the result is a $P$-genomic tableau, we 
say $\mu$ is {\bf admissible for $S$}; otherwise $\Delta_\mu(S)$ is not defined. 
Clearly, if $\Delta_\mu(S)$ is defined, it has content $\mu$. 

\begin{example}\label{ex:psi_type_B}
Let $S$ be the increasing tableau 
of Example~\ref{ex:phi_type_B}. Let $\eta = (2,1,1)$.
We compute $\Delta_\eta(S)$ in stages:
\[\begin{picture}(400,50)
\put(0,40){
\begin{ytableau}
*(lightgray)\blank & *(lightgray)\blank & *(lightgray)\blank & *(lightgray)\blank &  1_1 & 2_3\\
\none & *(lightgray)\blank & *(lightgray)\blank &  1_2 &  2_3\\
\none & \none &  3_4\\
\end{ytableau}}
\put(100,30){$\Longrightarrow$}
\put(240,30){$\Longrightarrow$}
\put(140,40){
\begin{ytableau}
*(lightgray)\blank & *(lightgray)\blank & *(lightgray)\blank & *(lightgray)\blank &  1' & 2\\
\none & *(lightgray)\blank & *(lightgray)\blank &  1 &  2\\
\none & \none &  3\\
\end{ytableau}}
\put(280,40){
\begin{ytableau}
*(lightgray)\blank & *(lightgray)\blank & *(lightgray)\blank & *(lightgray)\blank & *(red) 1' & *(green) 2\\
\none & *(lightgray)\blank & *(lightgray)\blank & *(SkyBlue) 1 & *(green) 2\\
\none & \none & *(Dandelion) 3\\ 
\end{ytableau}}
\end{picture}
\]
Observe that we obtain the genomic $P$-tableau $T$ of Example~\ref{ex:B.genomic.tableau}.

Compare this to the computation of $\Delta_{\theta}(S)$, where $\theta = (4)$: 
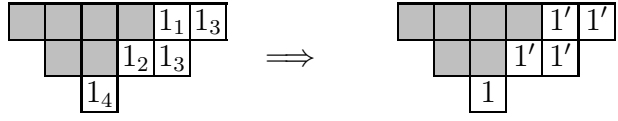
\[\begin{picture}(200,50)
\put(0,40){
\begin{ytableau}
*(lightgray)\blank & *(lightgray)\blank & *(lightgray)\blank & *(lightgray)\blank &  1_1 & 1_3\\
\none & *(lightgray)\blank & *(lightgray)\blank &  1_2 &  1_3\\
\none & \none &  1_4\\ 
\end{ytableau}}
\put(100,30){$\Longrightarrow$}
\put(150,40){\begin{ytableau}
*(lightgray)\blank & *(lightgray)\blank & *(lightgray)\blank & *(lightgray)\blank &  1' & 1'\\
\none & *(lightgray)\blank & *(lightgray)\blank &  1' &  1'\\
\none & \none &  1\\
\end{ytableau}}
\end{picture}
\]
Since the tableau obtained is not a $P$-tableau (it violates (P.3)), $\Delta_\theta(S)$ is undefined.
\qed
\end{example}

\begin{example}
Let $V$ be the increasing tableau $\ytableaushort{{*(lightgray)\blank} {*(lightgray)\blank} 1, \none 1 2}$ and let $\kappa = (2)$. Then in the construction of $\Delta_\kappa(V)$, we first obtain a valid $P$-tableau: \[\ytableaushort{ {*(lightgray)\blank} {*(lightgray)\blank} {1_1}, \none {1_1} {1_2} } \;\;\;\;\; \Longrightarrow \;\;\;\; \ytableaushort{{*(lightgray)\blank} {*(lightgray)\blank} {1'}, \none 1 1}.\]
However the putative genomic structure 
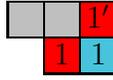
\[
\ytableaushort{{*(lightgray)\blank} {*(lightgray)\blank} {*(red) 1'}, \none {*(red) 1} {*(SkyBlue) 1}}
\]
is invalid, so $\Delta_\kappa(V)$ is undefined.
\qed
\end{example}

An increasing tableau $S$ is {\bf $\mu$-Pieri-filled} if $\mu$ is admissible for $S$ and $\Gamma(\Delta_\mu(S)) = S$. 

\begin{remark}
\label{rem:Pieri_filling}
It is easy to check that for $\mu = (\mu_1, \mu_2, \dots, \mu_h)$, an increasing tableau $S$ is $\mu$-Pieri-filled if and only if for each $k \leq h$, the entries of $S$ in $\mathcal{P}_k(\mu)$
form a \emph{Pieri filling} of a \emph{ribbon} in the sense of~\cite[$\mathsection$4]{Clifford.Thomas.Yong}. \qed
\end{remark}

\begin{lemma}\label{lem:shifted_bijection}
Let $T \in \Pgenomictableau{\nu/\lambda}{\mu}$. Then $\mu$ is admissible for $\Gamma(T)$ and $\Delta_\mu(\Gamma(T)) = T$. Hence $\Gamma(T)$ is $\mu$-Pieri-filled. 
\end{lemma}
\begin{proof}
The construction of $\Delta_{\mu}(\Gamma(T))$ is in stages. First we construct the underlying putative $P$-tableau structure for $\Delta_{\mu}(\Gamma(T))$. We will show that this is the same as the underlying $P$-tableau of $T$.
Consider a box ${\sf b}$ in $\nu/\lambda$. Suppose the box ${\sf b}$ in $T$ contains $\alpha \in \mathcal{D}$ (the color being irrelevant for now). Then it is clear that $\Delta_\mu(\Gamma(T))$ has $\beta \in {\sf b}$ with $|\beta| = |\alpha|$. The letter $\beta$ is primed if and only if there is $\gamma$ in box ${\sf c}$ southwest of ${\sf b}$ in $\Delta_\mu(\Gamma(T))$ with $|\gamma| = |\beta|$ and the entry of ${\sf c}$ in $\Gamma(T)$ strictly greater than the entry of ${\sf b}$ in $\Gamma(T)$. The entry of ${\sf c}$ in $\Gamma(T)$ is strictly greater than the entry of ${\sf b}$ in $\Gamma(T)$ exactly when ${\sf b} \prec {\sf c}$. By definition, this happens if and only if $\alpha$ is primed. Thus $\alpha = \beta$. Therefore $T$ and (the partially constructed tableau) $\Delta_\mu(\Gamma(T))$ have the same underlying $P$-tableau structure. 

In the next stage of constructing $\Delta_\mu(\Gamma(T))$, we attempt to partition the boxes into genes to produce a genomic $P$-tableau. By construction, $T$ and $\Delta_\mu(\Gamma(T))$ have the same partition of labels into genes; hence $\Delta_\mu(\Gamma(T))$ is defined and the first claim of the lemma holds. The second claim follows from the first by applying $\Gamma$.
\end{proof}

Let $\pieri{\nu / \lambda}{\mu} := \{ S : \text{$S$ is increasing of shape $\nu / \lambda$ and $\mu$-Pieri filled} \}$.

\begin{theorem}
$\Gamma : \Pgenomictableau{\nu / \lambda}{\mu} \to \pieri{\nu / \lambda}{\mu}$ and $\Delta_\mu : \pieri{\nu / \lambda}{\mu} \to \Pgenomictableau{\nu / \lambda}{\mu}$ are mutually inverse bijections.
\end{theorem}
\begin{proof}
Immediate by definition and Lemma~\ref{lem:shifted_bijection}.
\end{proof}

\subsection{Genomic $P$-Knuth equivalence}
Given a colored sequence $w$ of symbols from $\mathcal{D}$, write $\hat{w}$ for the sequence given by writing $w$ backwards, replacing each $k'$ with $k$ and each $k$ with $(k+1)'$ and preserving the colors (cf. the \emph{uncolored} definition of 
$\hat{w}$ after Example~\ref{ex:B.genomic.tableau}).
A {\bf genomic $P$-word} is a word $s$ of colored symbols from $\mathcal{D}$ such that in the concatenation $s \hat{s}$ all unprimed $i$'s of a fixed color are consecutive among the set of all unprimed $i$'s. Let ${\tt genomicseq}(T)$ denote the colored row reading word (right to left, and top to bottom) of a genomic $P$-tableau $T$, as for genomic tableaux in Section~\ref{sec:sequences}.

\begin{lemma}
Let $T$ be a genomic $P$-tableau. Then ${\tt genomicseq}(T)$ is a genomic $P$-word.
\end{lemma}
\begin{proof}
The follows from the fact that $T$ is a $P$-tableau together with the condition that the
boxes of each gene of family $i$ are consecutive in $\prec$-order.
\end{proof}

A {\bf genotype} of a genomic $P$-word $w$ is an uncolored subword given by choosing one letter of each color. A {\bf $P$-genotype} of the double sequence $w\hat{w}$ is a word of the form $x\hat{x}$ where $x$ is any genotype of $w$. We say $w\hat{w}$ is {\bf locally ballot} at the letter $\alpha$ if every $P$-genotype of $w\hat{w}$ that includes that $\alpha$ is locally ballot there. Finally $w\hat{w}$ is {\bf ballot} if every $P$-genotype of $w\hat{w}$ is ballot, equivalently if $w\hat{w}$ is locally ballot at each letter. In particular, the genomic $P$-tableau $T$ is ballot exactly when ${\tt genomicseq}(T)\widehat{{\tt genomicseq}(T)}$ is. 

\begin{example}
Let $T$ be the genomic $P$-tableau
$\begin{ytableau}
*(lightgray)\blank & *(lightgray)\blank & *(lightgray)\blank & *(lightgray)\blank & *(red) 1' & *(green) 2\\
\none & *(lightgray)\blank & *(lightgray)\blank & *(cyan) 1 & *(green) 2\\
\none & \none & *(Dandelion) 3\\
\end{ytableau}
$
of Example~\ref{ex:B.genomic.tableau}. Then 
\[{\tt genomicseq}(T)\widehat{{\tt genomicseq}(T)} = \textcolor{green}{2} \textcolor{red}{1'} \textcolor{green}{2} \textcolor{cyan}{1} \textcolor{Dandelion}{3} \textcolor{Dandelion}{4'} \textcolor{cyan}{2'} \textcolor{green}{3'} \textcolor{red}{1} \textcolor{green}{3'}.\] 
It has exactly two $P$-genotypes: 
\[\textcolor{green}{2} \textcolor{red}{1'} \textcolor{cyan}{1} \textcolor{Dandelion}{3} \textcolor{Dandelion}{4'} \textcolor{cyan}{2'} \textcolor{red}{1} \textcolor{green}{3'}
\text{\ and \ } \textcolor{red}{1'} \textcolor{green}{2} \textcolor{cyan}{1} \textcolor{Dandelion}{3} \textcolor{Dandelion}{4'} \textcolor{cyan}{2'} \textcolor{green}{3'} \textcolor{red}{1}.\] Neither $P$-genotype is ballot.
\qed
\end{example}

We define the equivalence relation $\equiv_{GP}$ of 
{\bf genomic $P$-Knuth equivalence} on genomic $P$-words as the transitive closure of the following relations:
\begin{align}
\mathbf{u}{\color{red}{\alpha \alpha}}\mathbf{v} &\equiv_{GP} \mathbf{u}{\color{red}{\alpha}}\mathbf{v},\tag{GP.1} \\ 
\mathbf{u}{\color{red}{\alpha}}{\color{SkyBlue}{\beta}}{\color{red}{\alpha}}\mathbf{v} &\equiv_{GP} \mathbf{u}{\color{SkyBlue}{\beta}}{\color{red}{\alpha}}{\color{cyan}{\beta}}\mathbf{v},\tag{GP.2} \\
\mathbf{u}{\color{SkyBlue}{\beta}}{\color{red}{\alpha}}{\color{green}{\gamma}}\mathbf{v} &\equiv_{GP} \mathbf{u}{\color{SkyBlue}{\beta}}{\color{green}{\gamma}}{\color{red}{\alpha}}\mathbf{v} \; \;\text{ if $\alpha \leq \beta < \gamma$ and $\beta=|\beta|$, or $\alpha < \beta \leq \gamma$ and $\beta=|\beta|'$,}\tag{GP.3} \\
\mathbf{u}{\color{red}{\alpha}}{\color{green}{\gamma}}{\color{SkyBlue}{\beta}}\mathbf{v} &\equiv_{GP} \mathbf{u}{\color{green}{\gamma}}{\color{red}{\alpha}}{\color{SkyBlue}{\beta}}\mathbf{v} \;\; \text{ if $\alpha \leq \beta < \gamma$ and $\beta=|\beta|'$, or $\alpha < \beta \leq \gamma$ and $\beta=|\beta|$,}\tag{GP.4} \\
\mathbf{u}{\color{red}{i}}{\color{SkyBlue}{j}} &\equiv_{GP} \mathbf{u}{\color{SkyBlue}{j^\dagger}}{\color{red}{i}}, \;\; \text{ where $j^\dagger = j'$ if $i=j$, and $j^\dagger = j$ otherwise,}\tag{GP.5}
\end{align}
where {\color{red}{red}}, {\color{SkyBlue}{blue}}, {\color{green}{green}} represent distinct colors. 

\begin{theorem}\label{prop:shifted_ballotness}
If $w_1 \equiv_{GP} w_2$, then $w_1\hat{w_1}$ is ballot if and only if $w_2\hat{w_2}$ is ballot.
\end{theorem}
\begin{proof}
Let $w$ be a genomic $P$-word. We need that (GP.1)--(GP.5) preserve ballotness of $w\hat{w}$. 

\noindent
(GP.1) and (GP.2):
These relations change $w$ without changing the set of genotypes of $w\hat{w}$. Hence they do not affect ballotness of the latter.

\noindent
(GP.3):
Suppose $w = \mathbf{u}{\color{SkyBlue}{\beta}}{\color{red}{\alpha}}{\color{green}{\gamma}}\mathbf{v}$ and that $w^* = \mathbf{u}{\color{SkyBlue}{\beta}}{\color{green}{\gamma}}{\color{red}{\alpha}}\mathbf{v}$ is obtained by (GP.3). 

(``$\Rightarrow$'' for (GP.3)): We assume
$w\hat{w}$ is ballot, and we must show  
$w^*\hat{w^*}$ is ballot, i.e., locally ballot (henceforth
abbreviated ``LB'') at each letter. 

{\sf (Case 1: $\alpha = \beta$)}: We have $i := |\alpha| = \alpha$ and $i < \gamma$. 

{\sf (Case 1.1: $\gamma = |\gamma|$)}: 
Let $k := \gamma$.
Then 
$$w\hat{w}=\mathbf{u}{\color{SkyBlue}{i}}{\color{red}{i}}{\color{green}{k}}\mathbf{v}\mathbf{\hat{v}}{\color{green}{(k+1)'}}{\color{red}{(i+1)'}}{\color{SkyBlue}{(i+1)'}}\mathbf{\hat{u}}$$ 
and
$$w^*\hat{w^*} = \mathbf{u}{\color{SkyBlue}{i}}{\color{green}{k}}{\color{red}{i}}\mathbf{v}\mathbf{\hat{v}}{\color{red}{(i+1)'}}{\color{green}{(k+1)'}}{\color{SkyBlue}{(i+1)'}}\mathbf{\hat{u}}.$$ 
It suffices to show that $w^*\hat{w^*}$ is LB at ${\color{green}{k}}$ and ${\color{green}{(k+1)'}}$. LBness at the latter is clear from the ballotness of $w\hat{w}$.

If $k > i+1$, then LBness at ${\color{green}{k}}$ is also clear from the ballotness of $w\hat{w}$. Hence assume $k = i+1$. The proof is now the same is for the corresponding case of Theorem~\ref{thm:knuth_preserves_ballot}.

{\sf (Case 1.2: $\gamma = |\gamma|'$)}: 
Let $k' := \gamma$.
Then $$w\hat{w}=\mathbf{u}{\color{SkyBlue}{i}}{\color{red}{i}}{\color{green}{k'}}\mathbf{v}\mathbf{\hat{v}}{\color{green}{k}}{\color{red}{(i+1)'}}{\color{SkyBlue}{(i+1)'}}\mathbf{\hat{u}}$$
and
$$w^*\hat{w^*} = \mathbf{u}{\color{SkyBlue}{i}}{\color{green}{k'}}{\color{red}{i}}\mathbf{v}\mathbf{\hat{v}}{\color{red}{(i+1)'}}{\color{green}{k}}{\color{SkyBlue}{(i+1)'}}\mathbf{\hat{u}}.$$ 
It suffices to show that $w^*\hat{w^*}$ is LB at ${\color{green}{k'}}$ and ${\color{green}{k}}$. LBness at the latter is clear from the ballotness of $w\hat{w}$. LBness at the former may be argued exactly as in {\sf Case 1.1}.

{\sf (Case 2: $\beta = \gamma$)}: We have $j' := |\beta|' = \beta$ and $\alpha < j'$. 

{\sf (Case 2.1: $\alpha = |\alpha|$)}: Let $i = \alpha$. 
Then 
$$w\hat{w}=\mathbf{u}{\color{SkyBlue}{j'}}{\color{red}{i}}{\color{green}{j'}}\mathbf{v}\mathbf{\hat{v}}{\color{green}{j}}{\color{red}{(i+1)'}}{\color{SkyBlue}{j}}\mathbf{\hat{u}}$$ 
and
$$w^*\hat{w^*} = \mathbf{u}{\color{SkyBlue}{j'}}{\color{green}{j'}}{\color{red}{i}}\mathbf{v}\mathbf{\hat{v}}{\color{red}{(i+1)'}}{\color{green}{j}}{\color{SkyBlue}{j}}\mathbf{\hat{u}}.$$ 
It suffices to show that $w^*\hat{w^*}$ is LB at ${\color{green}{j'}}$ and ${\color{green}{j}}$. 
That $w^*\hat{w^*}$ is LB at ${\color{green}{j'}}$ follows from the LBness of $w\hat{w}$ at ${\color{SkyBlue}{j'}}$. LBness at ${\color{green}{j}}$ is trivial.

{\sf (Case 2.2: $\alpha = |\alpha|'$)}: Let $i' = \alpha$.
Then $$w\hat{w}=\mathbf{u}{\color{SkyBlue}{j'}}{\color{red}{i'}}{\color{green}{j'}}\mathbf{v}\mathbf{\hat{v}}{\color{green}{j}}{\color{red}{i}}{\color{SkyBlue}{j}}\mathbf{\hat{u}}$$ 
and
$$w^*\hat{w^*} = \mathbf{u}{\color{SkyBlue}{j'}}{\color{green}{j'}}{\color{red}{i'}}\mathbf{v}\mathbf{\hat{v}}{\color{red}{i}}{\color{green}{j}}{\color{SkyBlue}{j}}\mathbf{\hat{u}}.$$ 
It suffices to check that $w^*\hat{w^*}$ is LB at ${\color{green}{j'}}$ and ${\color{green}{j}}$. 
This is clear from the ballotness of $w\hat{w}$.

{\sf (Case 3: $\alpha < \beta < \gamma$)}:

{\sf (Case 3.1: $\alpha = |\alpha|$)}: Let $i = \alpha$. If $\gamma > i+1$, ballotness is clear. Otherwise, by the assumptions of {\sf Case 3}, $\beta = (i+1)'$ and $\gamma = i+1$. 
So \[w\hat{w}=\mathbf{u}{\color{SkyBlue}{(i+1)'}}{\color{red}{i}}{\color{green}{(i+1)}}\mathbf{v}\mathbf{\hat{v}}{\color{green}{(i+2)'}}{\color{red}{(i+1)'}}{\color{SkyBlue}{(i+1)}}\mathbf{\hat{u}}\]
and
$$w^*\hat{w^*} = \mathbf{u}{\color{SkyBlue}{(i+1)'}}{\color{green}{(i+1)}}{\color{red}{i}}\mathbf{v}\mathbf{\hat{v}}{\color{red}{(i+1)'}}{\color{green}{(i+2)'}}{\color{SkyBlue}{(i+1)}}\mathbf{\hat{u}}.$$ 
It suffices to check that $w^*\hat{w^*}$ is LB at ${\color{green}{(i+1)}}$ and ${\color{green}{(i+2)'}}$.
The latter is clear from ballotness of $w\hat{w}$. The LBness at ${\color{green}{(i+1)}}$ follows from the LBness of $w\hat{w}$ at ${\color{SkyBlue}{(i+1)'}}$.

{\sf (Case 3.2: $\alpha = |\alpha|'$)}: Let $i' = \alpha$. If $\gamma > i+1$, ballotness is clear. Otherwise we have either $\gamma = (i+1)'$ or $\gamma = i+1$.

{\sf (Case 3.2.1: $\gamma = (i+1)'$)}: We have $\beta = i$.
Then $$w\hat{w}=\mathbf{u}{\color{SkyBlue}{i}}{\color{red}{i'}}{\color{green}{(i+1)'}}\mathbf{v}\mathbf{\hat{v}}{\color{green}{(i+1)}}{\color{red}{i}}{\color{SkyBlue}{(i+1)'}}\mathbf{\hat{u}}$$
and
$$w^*\hat{w^*} = \mathbf{u}{\color{SkyBlue}{i}}{\color{green}{(i+1)'}}{\color{red}{i'}}\mathbf{v}\mathbf{\hat{v}}{\color{red}{i}}{\color{green}{(i+1)}}{\color{SkyBlue}{(i+1)'}}\mathbf{\hat{u}}.$$ 
It suffices to check LBness at the two {\color{green} green} letters. These checks hold by the ballotness of $w\hat{w}$.

{\sf (Case 3.2.2: $\gamma = i+1$)}:
Here $$w\hat{w}=\mathbf{u}{\color{SkyBlue}{\beta}}{\color{red}{i'}}{\color{green}{(i+1)}}\mathbf{v}\mathbf{\hat{v}}{\color{green}{(i+2)'}}{\color{red}{i}}{\color{SkyBlue}{\hat{\beta}}}\mathbf{\hat{u}}$$
and
$$w^*\hat{w^*} = \mathbf{u}{\color{SkyBlue}{\beta}}{\color{green}{(i+1)}}{\color{red}{i'}}\mathbf{v}\mathbf{\hat{v}}{\color{red}{i}}{\color{green}{(i+2)'}}{\color{SkyBlue}{\hat{\beta}}}\mathbf{\hat{u}}.$$ 
It suffices to check LBness at the two {\color{green} green} letters. These checks are both direct from the ballotness of $w\hat{w}$.

(``$\Leftarrow$'' for (GP.3)): Conversely, assume
$w^*\hat{w^*}$ is ballot. We need to show that
 $w\hat{w}$ is ballot. As with the arguments for $\Rightarrow$, we 
need to establish LBness at each letter. In brief, it suffices to check this in each case below
at the {\color{green} green} letters. In each of these
situations, this is immediate from the assumption
$w^*\hat{w^*}$ is ballot.

{\sf (Case 1: $\alpha = \beta$)}: We have $i := |\alpha| = \alpha$ and $i < \gamma$. 

{\sf (Case 1.1: $\gamma = |\gamma|$)}: 
Let $k := \gamma$.
Then 
$$w^*\hat{w^*} = \mathbf{u}{\color{SkyBlue}{i}}{\color{green}{k}}{\color{red}{i}}\mathbf{v}\mathbf{\hat{v}}{\color{red}{(i+1)'}}{\color{green}{(k+1)'}}{\color{SkyBlue}{(i+1)'}}\mathbf{\hat{u}}$$ 
and
$$w\hat{w}=\mathbf{u}{\color{SkyBlue}{i}}{\color{red}{i}}{\color{green}{k}}\mathbf{v}\mathbf{\hat{v}}{\color{green}{(k+1)'}}{\color{red}{(i+1)'}}{\color{SkyBlue}{(i+1)'}}\mathbf{\hat{u}}.$$  

{\sf (Case 1.2: $\gamma = |\gamma|'$)}: 
Let $k' := \gamma$. 
Then 
\[w^*\hat{w^*} = \mathbf{u}{\color{SkyBlue}{i}}{\color{green}{k'}}{\color{red}{i}}\mathbf{v}\mathbf{\hat{v}}{\color{red}{(i+1)'}}{\color{green}{k}}{\color{SkyBlue}{(i+1)'}}\mathbf{\hat{u}}\] 
and
\[w\hat{w} = \mathbf{u}{\color{SkyBlue}{i}}{\color{red}{i}}{\color{green}{k'}}\mathbf{v}\mathbf{\hat{v}}{\color{green}{k}}{\color{red}{(i+1)'}}{\color{SkyBlue}{(i+1)'}}\mathbf{\hat{u}}.\]

{\sf (Case 2: $\beta = \gamma$)}: We have $j' := |\beta|' = \beta$ and $\alpha < j'$. 

{\sf (Case 2.1: $\alpha = |\alpha|$)}: Let $i = \alpha$. 
Then 
\[w^*\hat{w^*} = \mathbf{u}{\color{SkyBlue}{j'}}{\color{green}{j'}}{\color{red}{i}}\mathbf{v}\mathbf{\hat{v}}{\color{red}{(i+1)'}}{\color{green}{j}}{\color{SkyBlue}{j}}\mathbf{\hat{u}}\] 
and
\[w\hat{w} = \mathbf{u}{\color{SkyBlue}{j'}}{\color{red}{i}}{\color{green}{j'}}\mathbf{v}\mathbf{\hat{v}}{\color{green}{j}}{\color{red}{(i+1)'}}{\color{SkyBlue}{j}}\mathbf{\hat{u}}.\] 

{\sf (Case 2.2: $\alpha = |\alpha|'$)}: Let $i' = \alpha$.
Then 
\[w^*\hat{w^*} = \mathbf{u}{\color{SkyBlue}{j'}}{\color{green}{j'}}{\color{red}{i'}}\mathbf{v}\mathbf{\hat{v}}{\color{red}{i}}{\color{green}{j}}{\color{SkyBlue}{j}}\mathbf{\hat{u}}\]
and 
\[w\hat{w} = \mathbf{u}{\color{SkyBlue}{j'}}{\color{red}{i'}}{\color{green}{j'}}\mathbf{v}\mathbf{\hat{v}}{\color{green}{j}}{\color{red}{i}}{\color{SkyBlue}{j}}\mathbf{\hat{u}}.\] 

{\sf (Case 3: $\alpha < \beta < \gamma$)}:

{\sf (Case 3.1: $\alpha = |\alpha|$)}: Let $i = \alpha$. If $\gamma > i+1$, ballotness is clear. Otherwise, by the assumptions of {\sf Case 3}, we have $\beta = (i+1)'$ and $\gamma = i+1$. 
Thus 
\[w^*\hat{w^*} = \mathbf{u}{\color{SkyBlue}{(i+1)'}}{\color{green}{(i+1)}}{\color{red}{i}}\mathbf{v}\mathbf{\hat{v}}{\color{red}{(i+1)'}}{\color{green}{(i+2)'}}{\color{SkyBlue}{(i+1)}}\mathbf{\hat{u}}\]
and 
\[w\hat{w} = \mathbf{u}{\color{SkyBlue}{(i+1)'}}{\color{red}{i}}{\color{green}{(i+1)}}\mathbf{v}\mathbf{\hat{v}}{\color{green}{(i+2)'}}{\color{red}{(i+1)'}}{\color{SkyBlue}{(i+1)}}\mathbf{\hat{u}}.\]

{\sf (Case 3.2: $\alpha = |\alpha|'$)}: Let $i' = \alpha$. If $\gamma > i+1$, ballotness is clear. Otherwise we have either $\gamma = (i+1)'$ or $\gamma = i+1$.

{\sf (Case 3.2.1: $\gamma = (i+1)'$)}: We have $\beta = i$.
Then 
\[w^*\hat{w^*} = \mathbf{u}{\color{SkyBlue}{i}}{\color{green}{(i+1)'}}{\color{red}{i'}}\mathbf{v}\mathbf{\hat{v}}{\color{red}{i}}{\color{green}{(i+1)}}{\color{SkyBlue}{(i+1)'}}\mathbf{\hat{u}}\] 
and
\[w\hat{w} = \mathbf{u}{\color{SkyBlue}{i}}{\color{red}{i'}}{\color{green}{(i+1)'}}\mathbf{v}\mathbf{\hat{v}}{\color{green}{(i+1)}}{\color{red}{i}}{\color{SkyBlue}{(i+1)'}}\mathbf{\hat{u}}.\]

{\sf (Case 3.2.2: $\gamma = i+1$)}:
Here 
\[w^*\hat{w^*} = \mathbf{u}{\color{SkyBlue}{\beta}}{\color{green}{(i+1)}}{\color{red}{i'}}\mathbf{v}\mathbf{\hat{v}}{\color{red}{i}}{\color{green}{(i+2)'}}{\color{SkyBlue}{\hat{\beta}}}\mathbf{\hat{u}}\] 
and
\[w\hat{w} = \mathbf{u}{\color{SkyBlue}{\beta}}{\color{red}{i'}}{\color{green}{(i+1)}}\mathbf{v}\mathbf{\hat{v}}{\color{green}{(i+2)'}}{\color{red}{i}}{\color{SkyBlue}{\hat{\beta}}}\mathbf{\hat{u}}.\]

\noindent
(GP.4): This may be argued exactly as for (GP.3).

\noindent
(GP.5):
Suppose $w = \mathbf{u}{\color{red}{i}}{\color{SkyBlue}{j}}$ and that $w^* = \mathbf{u}{\color{SkyBlue}{j^\dagger}}{\color{red}{i}}$ is obtained by (GP.5). By symmetry, we may assume $i \leq j$. We must show $w\hat{w}$ is ballot if and only if $w^*\widehat{w^*}$ is.

{\sf (Case 1: $i < j$)}:
Then 
\[w\hat{w} = \mathbf{u}{\color{red}{i}}{\color{SkyBlue}{j}}{\color{SkyBlue}{(j+1)'}}{\color{red}{(i+1)'}}\mathbf{\hat{u}},\] while 
\[w^*\widehat{w^*} = \mathbf{u}{\color{SkyBlue}{j}}{\color{red}{i}}{\color{red}{(i+1)'}}{\color{SkyBlue}{(j+1)'}}\mathbf{\hat{u}}.\] 

Suppose $w\hat{w}$ is ballot.
It suffices to check LBness of $w^*\widehat{w^*}$ at the two {\color{SkyBlue} blue} letters. LBness at ${\color{SkyBlue}{(j+1)'}}$ is clear from the assumed ballotness of $w\hat{w}$. LBness at ${\color{SkyBlue}{j}}$, for $j = i+1$, follows from the LBness of $w\hat{w}$ at ${\color{red}{(i+1)'}}$ (when $j\neq i+1$, the claim is clear).

Conversely suppose $w^*\widehat{w^*}$ is ballot. It suffices to check LBness of $w\hat{w}$ at the two {\color{SkyBlue} blue} letters; this is immediate.

{\sf (Case 2: $i=j$)}:
Then 
\[w\hat{w} = \mathbf{u}{\color{red}{i}}{\color{SkyBlue}{i}}{\color{SkyBlue}{(i+1)'}}{\color{red}{(i+1)'}}\mathbf{\hat{u}},\] while 
\[w^*\widehat{w^*} = \mathbf{u}{\color{SkyBlue}{i'}}{\color{red}{i}}{\color{red}{(i+1)'}}{\color{SkyBlue}{i}}\mathbf{\hat{u}}.\] 
It is straightforward that $w\hat{w}$ is ballot if and only if $w^*\widehat{w^*}$ is.
\end{proof}

{\bf Weak $K$-Knuth equivalence} on words is the symmetric, transitive closure of these relations \cite[Definition~7.6]{Buch.Samuel}:
\begin{align*}
\mathbf{u}aa\mathbf{v} &\equiv_{wK} \mathbf{u}a\mathbf{v}, \\
\mathbf{u}aba\mathbf{v} &\equiv_{wK} \mathbf{u}bab\mathbf{v}, \\
\mathbf{u}bac\mathbf{v} &\equiv_{wK} \mathbf{u}bca\mathbf{v} \;\; \text{ if $a < b<c$},\\
\mathbf{u}acb\mathbf{v} &\equiv_{wK} \mathbf{u}cab\mathbf{v}  \;\; \text{ if $a < b<c$}, \\
\mathbf{u}ab &\equiv_{wK} \mathbf{u}ba.
\end{align*}

\begin{lemma}\label{lem:typeB_genomic tableau_knuth}
For genomic $P$-words $u, v$ we have $u \equiv_{GP} v$ if and only if $\Gamma(u) \equiv_{wK} \Gamma(v)$.
\end{lemma}
\begin{proof}
This follows from applying $\Delta_\mu$ to the generating relations for weak $K$-Knuth equivalence for Pieri-filled words.
\end{proof}

\subsection{Shifted jeu de taquin and the conclusion of the proof}

The definitions of genomic {\it jeu de taquin} and $K$-{\it jeu de taquin} for shifted tableaux are analogous to the unshifted case. For details of shifted $K$-{\it jeu de taquin}, see \cite{Clifford.Thomas.Yong}. We sketch the modifications necessary for shifted genomic {\it jeu de taquin} and give an illustrative example. For each gene $\GG$ of family $k$, define the
operator $\mathtt{switch}_{\GG}^{\bullet}$ as follows:
If ${\sf b}$ is a box of $\GG$ in the tableau $T$ with a 
neighbor containing a $\bullet$, replace the $k$ or $k' \in {\sf b}$ with $\bullet$ and remove it from $\GG$. If ${\sf c}$ is a box of $T$ containing a $\bullet$ and with a $\GG$ neighbor, ${\sf c}$ is a box of $\GG$ in $\mathtt{switch}_{\GG}^{\bullet}(T)$; ${\sf c}$ has entry $k$ in $\mathtt{switch}_{\GG}^{\bullet}(T)$ if either of its $\GG$ neighbors in $T$ have entry $k$ or if ${\sf c}$ lies on the main diagonal; otherwise ${\sf c}$ has entry $k'$ in $\mathtt{switch}_{\GG}^{\bullet}(T)$. The other boxes of $T$ are the same in $\mathtt{switch}_{\GG}^{\bullet}(T)$. 

Index the genes of $T$ as \[\GG_1 < \GG_2 < \dots < \GG_{|\mu|}\] 
according to the total order on genes from Section~\ref{sec:shifted_standardization}.  Then 
\[\jdt{T}{I}:=\mathtt{switch}_{\GG_{|\mu|}}^{\bullet} \circ \dots \circ \mathtt{switch}_{\GG_2}^{\bullet} \circ \mathtt{switch}_{\GG_1}^{\bullet}(T^{\bullet})\] 
with the $\bullet$'s deleted. (This algorithm reduces to the classical \emph{jeu de taquin} for semistandard $P$-tableaux in the case each gene contains only a single box.)

\begin{example} Suppose $T^\bullet$ is the genomic tableau
$\begin{ytableau}
*(lightgray)\blank & *(lightgray)\blank  & \bullet & *(red) 1' & *(SkyBlue) 1 \\
\none &  \bullet &*(red) 1' &*(Dandelion) 2\\
\none  & \none &*(green) 1
\end{ytableau}$. Then
\[\mathtt{switch}^\bullet_{{\color{red}1}}(T^\bullet) = 
\begin{ytableau}
*(lightgray)\blank & *(lightgray)\blank & *(red) 1' & \bullet & *(SkyBlue) 1 \\
\none & *(red) 1 & \bullet & *(Dandelion) 2\\
\none & \none  & *(green) 1
\end{ytableau}, \ \ \ 
\mathtt{switch}_{\color{green}1}^{\bullet} \circ \mathtt{switch}_{\color{red}1}^{\bullet}(T^\bullet) = \begin{ytableau}
*(lightgray)\blank & *(lightgray)\blank & *(red) 1'  & \bullet  & *(SkyBlue) 1\\
 \none & *(red) 1 & *(green) 1  &*(Dandelion) 2\\
\none & \none & \bullet 
\end{ytableau},\ \ \text{and }\]
\begin{align*}
\mathtt{switch}_{\color{Dandelion}2}^{\bullet} \circ  \mathtt{switch}_{\color{SkyBlue}1}^{\bullet} \circ \mathtt{switch}_{\color{green}1}^{\bullet} \circ \mathtt{switch}_{\color{red}1}^{\bullet}(T^\bullet) &=  \mathtt{switch}_{\color{SkyBlue}1}^{\bullet} \circ \mathtt{switch}_{\color{green}1}^{\bullet} \circ \mathtt{switch}_{\color{red}1}^{\bullet}(T^\bullet) \\ 
&= \begin{ytableau}
*(lightgray)\blank & *(lightgray)\blank & *(red) 1' & *(SkyBlue) 1 & \bullet\\
\none & *(red) 1 & *(green)1 & *(Dandelion) 2\\
\none & \none & \bullet 
\end{ytableau}.\ \
\end{align*}
\[
\text{So } \jdt{T}{I} = \begin{ytableau}
*(lightgray)\blank & *(lightgray)\blank & *(red) 1' & *(SkyBlue) 1 \\
\none & *(red) 1 & *(green)1 & *(Dandelion) 2\\
\none & \none 
\end{ytableau}.\] \qed
\end{example}

Using this shifted genomic {\it jeu de taquin}, one can obtain shifted versions of \emph{genomic infusion} and \emph{genomic Bender-Knuth involutions}, analogous to the discussion of Section~\ref{sec:genomic_Schur}. This leads to a definition of \emph{genomic $P$-Schur functions}, symmetric functions that deform the classical $P$-Schur functions just as the genomic Schur functions of Section~\ref{sec:genomic_Schur} deform the classical Schur functions. We do not pursue these ideas further here.

Let $S_\mu$ denote the {\bf row superstandard tableau} of shifted shape $\mu$ (that is, the tableau whose first row has entries $1,2,3,\ldots,\mu_1$, and whose
second row has entries $\mu_1+1,\mu_2+2,\ldots,\mu_1+\mu_2$ etc.). 

\begin{example}
For $\mu = (4,2)$, $S_\mu = \ytableaushort{1234,\none 56}$.
\qed
\end{example}

Let 
\[T_\mu := \Delta_\mu(S_\mu)\] 
be the unique genomic $P$-tableau whose underlying $P$-tableau is the highest weight tableau of shifted shape $\mu$.
We recall some results that we need.

\begin{theorem}[{\cite[Theorem 7.8]{Buch.Samuel}}]
\label{thm:shifted_BS}
Let $S$ be a shifted increasing tableaux.
Then $S$ rectifies to $S_\mu$ if and only if ${\tt seq}(S) \equiv_{wK} {\tt seq}(S_\mu)$. \qed
\end{theorem}

Let $\incr{\nu / \lambda}{\mu}:=\{\text{shifted increasing tableaux of shape $\nu / \lambda$ that rectify to $S_\mu$}\}.$

\begin{theorem}[{\cite[Theorem~1.2]{Clifford.Thomas.Yong}}]
\label{thm:CTY}
$b_{\lambda, \mu}^\nu = 
(-1)^{|\nu| - |\lambda| - |\mu|} \times \#\incr{\nu / \lambda}{\mu}$.
\qed
\end{theorem}

By Theorem~\ref{thm:CTY}, it is enough to biject
 $\incr{\nu / \lambda}{\mu}$ and $\Pballot{\nu / \lambda}{\mu}$. We claim that the maps $\Gamma$ and $\Delta_\mu$ give the desired bijections.
It follows from Remark~\ref{rem:Pieri_filling} and~\cite[Proof of Theorem~1.1]{Clifford.Thomas.Yong} that  $\Delta_\mu$ is well-defined on $\incr{\nu / \lambda}{\mu}$.

Let $S \in \incr{\nu / \lambda}{\mu}$. By Theorem~\ref{thm:shifted_BS}, 
\[{\tt seq}(S) \equiv_{wK} {\tt seq}(S_\mu).\] 
By Lemma~\ref{lem:typeB_genomic tableau_knuth}, 
\[{\tt genomicseq}(\Delta_\mu(S)) \equiv_{GP} {\tt genomicseq}(\Delta_\mu(S_\mu)) = {\tt genomicseq}(T_\mu).\] 
Note ${\tt genomicseq}(T_\mu)$ is ballot. Hence by Theorem~\ref{prop:shifted_ballotness}, ${\tt genomicseq}(\Delta_\mu(S))$ is ballot. Thus 
\[\Delta_\mu(S) \in \Pballot{\nu / \lambda}{\mu}.\]

Conversely, if $T \in  \Pballot{\nu / \lambda}{\mu}$, then its genomic rectification is also ballot by Theorem~\ref{thm:genomic tableau_knuth_equivalence_determines_plactic_class}. Hence its genomic rectification is $T_\mu$. Therefore 
\[\Gamma(T) \in \incr{\nu / \lambda}{\mu}.\] 
This completes the proof.\qed

\section*{Acknowledgments} 
We thank Hugh Thomas for many conversations which helped to make this project possible. OP was
supported by an NSF Graduate Research Fellowship, and Illinois Distinguished Fellowship from the University of Illinois, and NSF MCTP grant DMS 0838434. AY was supported by NSF grants and a Helen Corley Petit fellowship at UIUC.

\end{document}